\documentclass[final]{siamart1116}
\usepackage{amssymb}
\usepackage{eurosym}
\usepackage{lipsum}
\usepackage{amsfonts}
\usepackage{graphicx}
\usepackage{epstopdf}
\usepackage{algorithmic}
\usepackage{amsopn}

\ifpdf
  \DeclareGraphicsExtensions{.eps,.pdf,.png,.jpg}
\else
  \DeclareGraphicsExtensions{.eps}
\fi
\numberwithin{theorem}{section}
\newcommand{\TheTitle}{Restricted Dirichlet-to-Neumann data}
\newcommand{\TheAuthors}{M.V.Klibanov, J. Li and W. Zhang}

\headers{\TheTitle}{\TheAuthors}

\ifpdf
\hypersetup{
  pdftitle={\TheTitle},
  pdfauthor={\TheAuthors}
}
\fi
\begin{document}

\title{{Electrical Impedance Tomography with Restricted Dirichlet-to-Neumann
Map Data}\thanks{%
Submitted to the editors DATE. 
\funding{The wrok of MVK was supported by US Army Research Laboratory and US Army Research
Office grant W911NF-15-1-0233 and by the Office of Naval Research grant
N00014-15-1-2330.}}}
\author{ Michael V. Klibanov\thanks{%
Corresponding author. Department of Mathematics and Statistics, University
of North Carolina at Charlotte, Charlotte, NC 28223, USA (} \and Jingzhi Li 
\thanks{%
Department of Mathematics, Southern University of Science and Technology
(SUSTech), 1088 Xueyuan Boulevard, University Town of Shenzhen, Xili,
Nanshan, Shenzhen, Guangdong Province, P.R.China (} \and Wenlong Zhang%
\thanks{%
Department of Mathematics, Southern University of Science and Technology
(SUSTech), 1088 Xueyuan Boulevard, University Town of Shenzhen, Xili,
Nanshan, Shenzhen, Guangdong Province, P.R.China (} }
\maketitle

\begin{abstract}
We propose a new numerical method to reconstruct the isotropic electrical
conductivity from measured restricted Dirichlet-to-Neumann map data in
electrical impedance tomography (EIT) model. "Restricted
Dirichlet-to-Neumann (DtN) map data" means that the Dirichlet and Neumann
boundary data for EIT are generated by a point source running either along
an interval of a straight line or along a curve located outside of the
domain of interest. We \textquotedblleft convexify" the problem via
constructing a globally strictly convex Tikhonov-like functional using a
Carleman Weight Function. In particular, two new Carleman estimates are
established. Global convergenceto the correct solution of the gradient
projection method for this functional is proven. Numerical examples
demonstrate a good performance of this numerical procedure.
\end{abstract}




\section{Introduction}

\label{sect:1}

We develop in this paper a new globally convergent numerical method of the
reconstruction of the internal electrical conductivity in the inverse
problem of Electrical Impedance Tomography (EIT). The main part of the paper
is devoted to the theory of this method. Next, numerical examples are
presented. A general analytical concept of this method was originally
proposed in the work \cite{KlibJIIP} of the first author. However, it was
not sufficiently specified in \cite{KlibJIIP} for the EIT case. Unlike the
conventional case of the Dirichlet-to-Neumann map (DtN) boundary data, it
was proposed in \cite{KlibJIIP} to use the so-called \textquotedblleft
restricted DtN map data\emph{"} on the boundary. In the case of restricted
DtN data, the number $d=2,3$ of free variables in the data equals the number
of free variables in the unknown conductivity coefficient in the $\mathbb{R}%
^{d}$ case. We achieve this via truncation of a Fourier-like series. Note
that the conventional DtN requires $m=4$ of free variables in the data in
the 3D case.

Any Coefficient Inverse Problem (CIP) is highly nonlinear and ill-posed. As
a result, a conventional least squares Tikhonov functional for a CIP is
non-convex. The latter means that, as a rule, that functional has many local
minima and ravines, see, e.g. \cite{Scales} for a good numerical example.
Hence, to obtain a good approximation for the exact solution of a CIP, one
should start iterations of the minimization method for this functional in a
small neighborhood of the exact solution. We call this \emph{local
convergence}. However, it is often unclear how to practically obtain such a
good first guess.

Unlike the conventional case, we \textquotedblleft convexify" the problem.
More precisely, we construct a weighted Tikhonov-like functional with the
Carleman Weight Function (CWF) in it. The CWF is the function which is
involved in the Carleman estimate for the Laplace operator. The presence of
the CWF ensures the strict convexity of this functional on any \emph{a priori%
} chosen ball of an arbitrary radius $R>0$ in an appropriate Hilbert space.
The latter guarantees the global convergence of the gradient projection
method of the optimization of this functional to the exact solution of the
original inverse EIT problem. We call a numerical method for a CIP \emph{%
globally convergent} if there is a theorem, which guarantees that this
method delivers at least one point in a sufficiently small neighborhood of
the exact solution of that CIP without any advanced knowledge of this
neighborhood. The size of this neighborhood should depend on measurement and
approximation errors. The numerical method of this paper converges globally.

Electrical impedance tomography (EIT) is a non invasive and diffusive
imaging method to recover the electrical conductivity distribution inside an
object of interest by using the DtN map on the boundary. This modality is
safe, portable and also has many clinical imaging applications \cite%
{Holder:2005}. There is a vast number of research papers discussing EIT. It
has been analytically proven that the interior electrical conducting is
uniquely determined by the Dirichlet-to-Neumann map on the boundary \cite%
{Calderon:1980,Nachman,SU}. However, the EIT inverse problem is essentially
ill-posed compared with other imaging modalities in practice \cite{borcea},
since the DtN data on the boundary is not that sensitive to the conductivity
change inside the domain of interest.

In the past three decades, there were numerous studies on the EIT imaging
method with quite many publications. Since this paper is not a survey of
EIT, we now provide a far incomplete list of references on this topic. The
recovery of small inclusions from boundary measurements is discussed in \cite%
{anomaly1, anomaly4}. Hybrid conductivity imaging methods are presented in 
\cite{Wenlong2017,hybrid8,hybrid9}. The multi-frequency EIT imaging methods
are discussed in \cite{AmmariWenlong2016,SeoLeeKim:2008}. In particular, 
\cite{AmmariWenlong2016} also shows that the frequency difference method can
eliminate the modeling errors. Both the finite element method and the
adaptive finite element method are also applied to recover the internal
conductivity \cite{BangtiFEM2016,BangtiFEM2014}. The imaging algorithms
based on the sparsity reconstruction are considered in \cite%
{AmmariWenlong2016,Bangtisparse2012}. In \cite{Harrah} a globally convergent
method for shape reconstruction in EIT is proposed. Siltanen and Mueller
have done a lot work on the EIT inverse problem, including d-bar method,
diction reconstruction method, recovering boundary shape and imaging the
anisotropic electrical conductivity \cite%
{Mueller2016,Mueller2014,Siltanen2016,Siltanen2014}. Hyvonen, P$\ddot{a}$iv$%
\ddot{a}$rinta and Tamminen also offer in their recent paper a new way to
solve EIT problem \cite{Paivarinta2018}.

In a typical EIT experiment, constant electrical currents are applied to the
electrodes on the boundary of the object to image. Then the electrical
potentials are measured on the boundary. This gives the DtN map data. The
EIT problem is to recover the internal electric conductivity from these DtN
measurements. This problem is essentially ill-posed.

Unlike the DtN, by our definition, restricted DtN data means that the
Dirichlet and Neumann boundary data for the EIT problem are generated by a
point source running either along an interval of a straight line or along a
curve located outside of the domain to be imaged. Moreover, the restricted
DtN data can be given either on the whole or on a part of the boundary.

The key element of our method consists in the construction of a weighted
Tikhonov-like functional which is strictly convex on any a priori chosen
ball of an arbitrary radius $R>0$ in an appropriate Hilbert space. In other
words, we \textquotedblleft convexify" the problem. The main ingredient of
that Tikhonov-like functional is the presence of the CWF in it. If the exact
solution belongs to that ball (as it should be assumed in the framework of
the regularization theory \cite{T}), then convergence of the gradient
projection method to the exact solution is guaranteed if starting from an
arbitrary point of this ball. Hence, this is \emph{global convergence}. On
the other hand, recall that convergence of a gradient-like method to the
exact solution for a non-convex functional might be guaranteed only if its
starting point is located in a small neighborhood of this solution.

Carleman estimates were introduced in the field of Coefficient Inverse
Problems (CIPs) in the work \cite{BukhKlib}. There are many works devoted to
the method of \cite{BukhKlib}, see, e.g. the survey \cite{Ksurvey}, the most
recent book \cite{BY} and the references cited therein. The goal of the
authors of \cite{BukhKlib} was to prove global uniqueness and stability
results for CIPs. Later, however, it became clear that Carleman estimates
can also be applied to numerical methods for some ill-posed problems for
PDEs. First, CWFs can be applied to convexify CIPs, see \cite%
{BKconv,KlibSIMA97,KlibKam} for the theory and \cite%
{Baud,KlibThanh,KlibKol1,KlibKol2} for both the theory and numerical
results. Second, CWFs can be applied to prove convergence of the so-called
quasi-reversibility method for ill-posed Cauchy problems for linear PDEs 
\cite{KQR}. Third, CWFs are applicable for the convexification of some
ill-posed Cauchy problems for quasilinear PDEs, see \cite{KlibCauchy} for
the theory and \cite{BakKlib,KlibYag} for both the theory and numerical
results.

However, in the above cited works on the convexification for CIPs only the
case of a single location of the source was considered for either time
dependent or frequency dependent data. Unlike the above cited publications,
in \cite{KlibJIIP}\ a significantly new convexification method was proposed.
This was done for the case when the boundary data for a CIP are generated by
a point source which is running along an interval of a straight line. The
resulting boundary data form the above mentioned restricted DtN. In this
work we specify the idea of \cite{KlibJIIP} for the case of an inverse
problem for EIT with the restricted DtN data.

To minimize the above mentioned weighted Tikhonov-like functional, we
propose a multi-level method, which is somewhat similar with the adaptivity
method, see, e.g. \cite{BK} for a detailed theory of the adaptivity.
However, we do not extend to our case the theory of the adaptivity presented
in \cite{BK}, i.e. we restrict our attention only to the numerical aspect of
the adaptivity. Thus, we minimize that functional on a coarse mesh first and
use the solution achieved on the coarse mesh (first level) as the starting
point for a finer mesh (second level). We repeat this process until we get a
solution on $K_{th}$ level. We have found that we get a rough image on the
coarse mesh (e.g. support, shape) of the internal conductivity much faster
than on a finer mesh, while on the finer mesh with the starting point from
the solution on the coarse mesh, the solution is corrected in details (e.g.,
amplitude and shape).

\section{EIT with restricted Dirichlet-to-Neumann (DtN) data}

\label{sect:2}

All functions below are real valued ones. The same is about functional
spaces, including Hilbert spaces.

\subsection{Model}

\label{sect:2.1}

In this section, we formulate the restricted DtN for the inverse EIT
problem. First, we recall the traditional DtN for EIT. Let $\Omega $ be an
open bounded domain in $\mathbb{R}^{d}$ ($d=2,3$) to be imaged with a smooth
boundary $\partial \Omega $. The EIT forward problem is formulated as: For
any given input current%
\[
g_{1}\in L_{0}^{2}(\partial \Omega ):=\{g\in L^{2}(\Omega ):\int_{\partial
\Omega }g\,ds=0\} 
\]%
and the conductivity distribution $\sigma (x)$, find the function $u(x)\in
H^{1}(\Omega )$ such that 
\begin{equation}
\left\{ 
\begin{aligned} -\nabla\cdot(\sigma(x)\nabla u(x)) & = 0\quad \mbox{
in }\Omega,\\ \sigma(x)\frac{\partial u}{\partial\nu}& = g_1(x)\quad \mbox{on }\partial\Omega,\\ \int_{\partial\Omega}u(x)\,ds& =0, \end{aligned}%
\right.  \label{eqn:eit}
\end{equation}%
where $\nu $ is the outward unit normal vector on $\partial \Omega $. Denote 
$g_{0}(x)=u|_{\partial \Omega }$. Then the inverse EIT problem is to recover
the internal conductivity function $\sigma (x)$ from the DtN map $\Lambda
:g_{0}\rightarrow g_{1}$.

In this paper, we consider the EIT problem with the source outside the
domain of interest and the restricted DtN data measured on the boundary of
the domain of interest, as described below.

To avoid working with singularities and also to simplify the presentation,
we model the point source here by a $\delta -$like function instead of the $%
\delta -$function. Let $\varepsilon >0$ be a sufficiently small number. Let
the source function $f(x)$ be such that%
\begin{equation}
f(x)\in C^{\infty }(\mathbb{R}^{n}),f(0)\neq 0,f(x)\geq 0,\forall x\in 
\mathbb{R}^{d},\text{ }f(x)=0\text{ for }|x|>\varepsilon .  \label{2.100}
\end{equation}%
Let $G\subset \mathbb{R}^{n}$ be a bounded domain with its boundary $%
\partial G\in C^{1}$, $\Omega \subset G$ and $\partial \Omega \cap \partial
G=\varnothing $. Let $\overline{x}\in \mathbb{R}^{d-1}$ be a fixed point.
For $s\in \lbrack 0,1]$ denote $x_{s}=\left( x_{1s},\overline{x}\right) $
the position of the point source. Let $I=\{x_{s}=\left( x_{1s},\overline{x}%
\right) :s\in \lbrack 0,1]\}$ be the interval of the straight line $\left\{
x=\left( x_{1},\overline{x}\right) ,x_{1}\in \mathbb{R}\right\} .$ Let $%
I_{\varepsilon }=\{x\in \mathbb{R}^{d}:dist(x,I)<\varepsilon \}$, where $%
dist(x,I)$ is the Hausdorff distance between the point $x$ and the set $I$.
We also assume that $I_{\varepsilon }\subset (G\setminus \overline{\Omega }%
), $ which means that the support of the source function is outside of the
domain $\Omega $.

Let the function 
\begin{equation}
\sigma \left( x\right) \in C^{2+\alpha }(\overline{G}),\sigma \left(
x\right) =1\text{ for }x\in G\setminus {\Omega }\text{ and }\sigma \left(
x\right) \geq \sigma _{0}=const.>0.  \label{2.101}
\end{equation}%
Here $\alpha =const.\in \left( 0,1\right) $ and $C^{k+\alpha }(\overline{G})$
be the H\"{o}lder space, where $k\geq 0$ is an integer. Assume first that $%
\sigma (x)$ is known. For each source position $x_{s}\in I$ we define the
forward boundary value problem for EIT as the problem of finding the
function $u(x,s)$ such that 
\begin{equation}
\left\{ 
\begin{array}{lll}
\nabla \cdot (\sigma (x)\nabla u(x,s)) & =-f(x-x_{s}), & \qquad x\in
G,\forall x_{s}\in \overline{I} \\ 
u(x,s)|_{x\in \partial G} & =0, & \qquad \forall x_{s}\in \overline{I}.%
\end{array}%
\right.  \label{eqn:eitres}
\end{equation}%
It is well known that for each $x_{s}\in I$ the problem (\ref{eqn:eitres})
has unique solution 
\begin{equation}
u(x,s)\in C^{3+\alpha }(\overline{G}),\forall x_{s}\in \overline{I},
\label{2.106}
\end{equation}%
see, e.g. \cite{Gilbarg}. We measure both Dirichlet and Neumann boundary
conditions of the function $u$ on a part $\Gamma \subseteq \partial \Omega $
of the boundary $\partial \Omega $,

\begin{equation}
u(x,s)|_{x\in \Gamma ,x_{s}\in \overline{I}}=g_{0}(x,s)\text{ and }\partial
_{\nu }u(x,s)|_{x\in \Gamma ,x_{s}\in \overline{I}}=g_{1}(x,s).
\label{eitres-bd}
\end{equation}%
We call the Dirichlet and Neumann boundary data (\ref{eitres-bd})\emph{\ }
\textquotedblleft restricted DtN data".

If the coefficient $\sigma (x)$ is known, then, having the solution of the
forward problem (\ref{eqn:eitres}), one can easily compute functions $%
g_{0}(x,s)$ and $g_{1}(x,s).$ \ Suppose now that the function $\sigma (x)$
is unknown. Then we arrive at the following inverse problem:

\textbf{Coefficient Inverse Problem (CIP). }\emph{Assume that the function }$%
\sigma \left( x\right) $\emph{\ is unknown for }$x\in \Omega $\emph{\ and
also that conditions (\ref{2.100}), (\ref{2.101}) hold. Also, assume that
functions }$g_{0}(x,s)$\emph{\ and }$g_{1}(x,s)$\emph{\ in (\ref{eitres-bd})
are known for all }$x\in \Gamma ,x_{s}\in \overline{I}.$\emph{\ Determine
the function} $\sigma \left( x\right) .$

Note that in this CIP the number $d$ of free variables in the data equals
the number of free variables in the unknown coefficient.

\subsection{An equivalent problem}

\label{sect:2.2}

In this subsection, we transform the above CIP to an inverse problem for a
quasilinear PDE. First, introduce the well known change of variables 
\begin{equation}
u_{1}=\sqrt{\sigma }u,  \label{2.102}
\end{equation}%
where $u(x,s)$ is the solution of problem (\ref{eqn:eitres}). Then 
\begin{equation}
\left\{ 
\begin{array}{lll}
\Delta u_{1}(x,s)+a_{0}(x)u_{1}(x,s) & =-f(x-x_{s}), & \qquad \forall
x_{s}\in \overline{I}, \\ 
u_{1}(x,s)|_{x\in \partial G} & =0, & \qquad \forall x_{s}\in \overline{I},%
\end{array}%
\right.  \label{eqn:eitu1}
\end{equation}%
where 
\begin{equation}
{a_{0}(x)=-\frac{\Delta \left( \sqrt{\sigma (x)}\right) }{\sqrt{\sigma (x)}}}%
.  \label{2.103}
\end{equation}%
Recalling that $\sigma =1$ on $\partial \Omega ,$ we obtain from (\ref%
{eitres-bd}) 
\begin{equation}
u_{1}(x,s)|_{x\in \Gamma ,s\in \lbrack 0,1]}=g_{0}(x,s)\text{ and }\partial
_{n}u_{1}(x,s)|_{x\in \Gamma ,s\in \lbrack 0,1]}=g_{1}(x,s).  \label{2.104}
\end{equation}%
If we would recover the function $a_{0}(x)$ for $x\in \Omega $ from
conditions (\ref{eqn:eitu1}), (\ref{2.104}), then, assuming that $0$ is not
an eigenvalue of the elliptic operator $\Delta +a_{0}\left( x\right) $ with
the Dirichlet boundary condition either on $\partial \Omega $ or on $%
\partial G,$ we would recover the function $\sigma \left( x\right) $ via
solving elliptic equation (\ref{2.103}) either in the domain $\Omega $ with
the Dirichlet boundary condition $\sigma \mid _{\partial \Omega }=1,$ or in
the domain $G$ with the Dirichlet boundary condition $\sigma \mid _{\partial
G}=1.$ Hence, we focus below on the recovery of the function $a_{0}(x)$ for $%
x\in \Omega $ from conditions (\ref{eqn:eitu1}), (\ref{2.104}).

It follows from (\ref{2.100}), (\ref{eqn:eitres}), (\ref{2.102}) and the
maximum principle for elliptic equations \cite{Gilbarg} that $u_{1}(x,s)>0$
for all $x\in \overline{\Omega }$ and all $s\in \left[ 0,1\right] .$ Hence,
we can consider the function $v(x,s),$ 
\begin{equation}
v(x,s)=\ln u_{1}(x,s).  \label{2.105}
\end{equation}%
Then $u_{1}(x,s)=e^{v(x,s)}$ and (\ref{eqn:eitu1}) implies that 
\begin{equation}
\Delta v(x,s)+\left( \nabla v(x,s)\right) ^{2}=-a_{0}(x),\qquad x\in \Omega
,\forall s\in \lbrack 0,1].  \label{eqn:eitv}
\end{equation}%
Here we use (\ref{2.100}) and the fact that $I_{\varepsilon }\subset
(G\setminus \overline{\Omega })$. In addition, using (\ref{2.104}), we
obtain 
\begin{equation}
v(x,s)|_{x\in \Gamma ,s\in \lbrack 0,1]}=\tilde{g}_{0}(x,s)\text{ and }%
\partial _{\nu }v(x,s)|_{x\in \Gamma ,s\in \lbrack 0,1]}=\tilde{g}_{1}(x,s),
\label{eitv-bd}
\end{equation}%
where 
\begin{equation}
\tilde{g}_{0}(x,s)=\ln g_{0}(x,s)\text{ and }\tilde{g}_{1}(x,s)=\frac{%
g_{1}(x,s)}{g_{0}(x,s)}.  \nonumber
\end{equation}

Differentiating equation \eqref{eqn:eitv} with respect to $s$ and noting
that the function $a_{0}(x)$ is independent on $s$, we obtain 
\begin{equation}
\Delta v_{s}+2\nabla v_{s}\cdot \nabla v=0,\qquad x\in \Omega ,\forall s\in
\lbrack 0,1].  \label{eqn:eitvdiff}
\end{equation}%
Now the above CIP is reduced to the following problem:

\textbf{Reduced Problem}. \emph{Recover the function }$v\left( x,s\right) $%
\emph{\ from equations \eqref{eqn:eitvdiff}, given the boundary measurements 
}$\tilde{g}_{0}(x,s)$\emph{\ and }$\tilde{g}_{1}(x,s)$\emph{\ in %
\eqref{eitv-bd}.}

If the function $v\left( x,s\right) $ is approximated, then the approximate
coefficient $a_{0}(x)$ can be found via (\ref{eqn:eitv}). Thus, our focus
below is on the solution of Reduced Problem.

\section{Cauchy problem for a system of coupled quasilinear elliptic
equations}

\label{sect:3}

To solve Reduced Problem, we obtain in this section a system of coupled
quasilinear elliptic equations.

\subsection{A special orthonormal basis in $L^{2}\left( 0,1\right) $}

\label{sect:3.1}

Let $\left( ,\right) $ denotes the scalar product in $L^{2}\left( 0,1\right)
.$ We need to construct such an orthonormal basis in the space $L^{2}\left(
0,1\right) $ of real valued functions $\left\{ \psi _{n}\left( s\right)
\right\} _{n=0}^{\infty }$ that the following two conditions are met:

\begin{enumerate}
\item $\psi _{n}\in C^{1}\left[ 0,1\right] ,$ $\forall n=0,1,...$

\item Let $a_{mn}=\left( \psi _{n}^{\prime },\psi _{m}\right) .$ Then the
matrix $M_{k}=\left( a_{mn}\right) _{m,n=0}^{k-1}$ should be invertible for
any $k=1,2,...$
\end{enumerate}

Neither the basis of any type of classical orthonormal polynomials nor the
basis of trigonometric functions $\left\{ \sin \left( 2\pi ns\right) ,\cos
\left( 2\pi ns\right) \right\} _{n=0}^{\infty }$ do not satisfy the second
condition. Indeed, in either of these cases all elements of the first raw of
the matrix $M_{k}$ would be equal to zero. The required basis was
constructed in \cite{KlibJIIP}. We now briefly describe this construction
for the convenience of the reader.

Consider the set of functions $\left\{ s^{n}e^{s}\right\} _{n=0}^{\infty }.$
This set is complete in $L^{2}\left( 0,1\right) .$ We orthonormalize it
using the classical Gram-Schmidt orthonormalization procedure. We start from 
$n=0$, then take $n=1$, etc. Then we obtain the orthonormal basis $\left\{
\psi _{n}\left( s\right) \right\} _{n=0}^{\infty }$ in $L^{2}\left(
0,1\right) .$ Each function $\psi _{n}\left( s\right) $ has the form%
\begin{equation}
\psi _{n}\left( s\right) =P_{n}\left( s\right) e^{s},  \label{3.1}
\end{equation}%
where $P_{n}\left( s\right) $ is the polynomial of the degree $n$. Hence,
one can say that these polynomials are orthogonal to each other in the
weighted $L^{2}\left( 0,1\right) $ space with the weight function $e^{2s}.$
Lemma 3.1 ensures that the above property number 2 holds for functions $\psi
_{n}\left( s\right) $.

\textbf{Lemma 3.1} \cite{KlibJIIP}\textbf{.} \emph{We have} 
\begin{equation}
a_{mn}=\left[ \psi _{n}^{\prime },\psi _{m}\right] =\left\{ 
\begin{array}{c}
1\text{ if }n=m, \\ 
0\text{ if }n<m.%
\end{array}%
\right.  \label{3.2}
\end{equation}%
\emph{For an integer }$k\geq 1$\emph{\ consider the }$k\times k$\emph{\
matrix }$M_{k}=\left( a_{mn}\right) _{\left( m,n\right) =\left( 0,0\right)
}^{\left( k-1,k-1\right) }.$\emph{\ Then (\ref{3.2}) implies that }$M_{k}$%
\emph{\ is an upper diagonal matrix and} $\det \left( M_{k}\right) =1.$ 
\emph{Thus, the inverse matrix }$M_{k}^{-1}$\emph{\ exists. }

\subsection{Cauchy problem for a system of coupled quasilinear elliptic
equations}

\label{sect:3.2}

Fix an integer $N\geq 1.$ Denote $\Psi \left( N\right) =\left\{ \psi
_{n}\left( s\right) \right\} _{n=0}^{N-1}.$ We assume that the function $%
v\left( x,s\right) $ in (\ref{2.105}) can be represented via the truncated
Fourier-like series with respect to the orthonormal basis of functions $\psi
_{n}\left( s\right) $ in (\ref{3.1}),%
\begin{equation}
v\left( x,s\right) =\mathop{\displaystyle \sum }\limits_{n=0}^{N-1}v_{n}%
\left( x\right) \psi _{n}\left( s\right) ,\text{ }x\in \Omega ,\forall s\in
\lbrack 0,1].  \label{3.3}
\end{equation}%
Then the derivative $v_{s}\left( x,s\right) $ is%
\begin{equation}
v_{s}\left( x,s\right) =\mathop{\displaystyle \sum }\limits_{n=0}^{N%
\_1}v_{n}\left( x\right) \psi _{n}^{\prime }\left( s\right) ,\text{ }x\in
\Omega ,\forall s\in \lbrack 0,1].  \label{3.4}
\end{equation}%
Note that functions $v_{n}\left( x\right) $ are unknown and should be
determined. By (\ref{2.106}) and (\ref{2.102}) it is reasonable to assume
that functions $v_{n}\left( x\right) $ are such that 
\begin{equation}
v_{n}\in C^{3}\left( \overline{\Omega }\right) ,n=0,...,N-1.  \label{3.40}
\end{equation}%
It is likely that (\ref{3.40}) can be proven using the classical theory of
elliptic PDEs \textbf{\ }\cite{Gilbarg}. However, we are not doing this here
for brevity.

Substituting (\ref{3.3}) and (\ref{3.4}) in (\ref{eqn:eitvdiff}), we obtain%
\begin{equation}
\mathop{\displaystyle \sum }\limits_{n=0}^{N-1}\Delta v_{n}\left( x\right)
\psi _{n}^{\prime }\left( s\right) +\mathop{\displaystyle \sum }%
\limits_{n,k=0}^{N-1}\nabla v_{n}\left( x\right) \nabla v_{k}\left( x\right)
\psi _{n}^{\prime }\left( s\right) \psi _{k}\left( s\right) =0,\text{ }x\in
\Omega ,\forall s\in \lbrack 0,1].  \label{3.5}
\end{equation}%
Consider the vector function of unknown coefficient $v_{n}\left( x\right) $
in the expansion (\ref{3.3}),%
\begin{equation}
V\left( x\right) =\left( v_{0}\left( x\right) ,...,v_{N-1}\left( x\right)
\right) ^{T}.  \label{3.6}
\end{equation}%
For $m=0,...,N-1$ multiply both sides of (\ref{3.5}) by the function $\psi
_{m}\left( s\right) $ and then integrate with respect to $s\in \left(
0,1\right) .$ Using (\ref{3.40}) and (\ref{3.6}), we obtain%
\begin{equation}
M_{N}\Delta V-\widetilde{F}\left( \nabla V\right) =0,\text{ }x\in \Omega
,V\in C^{3}\left( \overline{\Omega }\right) ,  \label{3.7}
\end{equation}%
where the $N-$dimensional vector function $\widetilde{F}$ is quadratic with
respect to the first derivatives $\partial _{x_{j}}v_{k}\left( x\right)
,j=1,..,d;k=0,...,N-1.$ Multiplying both sides of (\ref{3.7}) by the inverse
matrix $M_{N}^{-1}$ (Lemma 3.1), we obtain a system of coupled quasilinear
elliptic equations,%
\begin{equation}
\Delta V-F\left( \nabla V\right) =0,x\in \Omega ,V\in C^{3}\left( \overline{%
\Omega }\right) ,  \label{3.8}
\end{equation}
\begin{equation}
F\left( \nabla V\right) =M_{N}^{-1}\widetilde{F}\left( \nabla V\right) .
\label{3.80}
\end{equation}%
Since the vector function $\widetilde{F}$ is quadratic with respect to the
first derivatives $\partial _{x_{j}}v_{k}\left( x\right) ,$ then (\ref{3.80}%
) implies that the vector function $F$ is also quadratic. In addition, using
(\ref{eitv-bd}), we obtain Cauchy data for the vector function $V\left(
x\right) $ on $\Gamma ,$%
\begin{equation}
V\left( x\right) \mid _{\Gamma }=p_{0}\left( x\right) ,\partial _{\nu
}V\left( x\right) \mid _{\Gamma }=p_{1}\left( x\right) .  \label{3.9}
\end{equation}

If we would solve the Cauchy problem (\ref{3.8}), (\ref{3.9}), then we would
find coefficients $v_{n}\left( x\right) $ in (\ref{3.3}). Next, we would
substitute (\ref{3.3}) in (\ref{eqn:eitv}) and obtain the following
approximate formula for the function $a_{0}\left( x\right) :$%
\begin{equation}
a_{0}\left( x\right) =-\mathop{\displaystyle \sum }\limits_{n=0}^{N-1}\Delta
v_{n}\left( x\right) \psi _{n}\left( s\right) +\left( \mathop{\displaystyle
\sum }\limits_{n=0}^{N-1}\nabla v_{n}\left( x\right) \psi _{n}\left(
s\right) \right) ^{2},\text{ }x\in \Omega ,s\in \left( 0,1\right) .
\label{3.90}
\end{equation}
As to the value of the parameter $s$ for which the function $a_{0}\left(
x\right) $ should be calculated in (\ref{3.90}), it should be chosen
numerically, similarly with \cite{KlibThanh,KlibKol1,KlibKol2}. Hence, we
develop below a numerical method for solving problem (\ref{3.8}), (\ref{3.9}%
).

\subsection{Two new Carleman estimates}

\label{sect:3.3}

Since in our numerical examples the domain $\Omega \subset \mathbb{R}^{2}$
is a disk, we prove in this subsection a new Carleman estimate for the
Laplace operator, which is specifically used for the disk in the 2D case and
for the ball in the 3D case. We work with the case when $\Gamma =\partial
\Omega $ since this is done in our numerical experiments. In principle,
Carleman estimates are known for this kind of domains, see, e.g. \cite%
{KlibSIMA97}. However, the CWF in \cite{KlibSIMA97} has a rather complicated
form and changes too rapidly. On the other hand, the previous numerical
experience of the first author with the convexification for CIPs \cite%
{KlibThanh,KlibKol1,KlibKol2} tells us that one should use a CWF of the
simplest possible form, also, see, e.g. \cite{Baud} for a similar statement.
This is the reason of presenting here the Carleman estimate with a simple
CWF which was not used before.

\subsubsection{The 3D case}

\label{sect:3.3.1}

We derive in this section a new Carleman estimate for the 3D case when the
domain $\Omega $ is a ball of the radius $\rho $,%
\begin{equation}
\Omega =\left\{ x\in \mathbb{R}^{3}:\left\vert x\right\vert <\rho \right\} .
\label{3.900}
\end{equation}%
Let $\mu \in \left( 0,\rho \right) $ be a number. Define the domain $\Omega
_{\mu }$ as%
\begin{equation}
\Omega _{\mu }=\left\{ x\in \mathbb{R}^{3}:\mu <\left\vert x\right\vert
<\rho \right\} \subset \Omega .  \label{3.901}
\end{equation}%
Consider spherical coordinates%
\[
r=\left\vert x\right\vert \in \left( \mu ,\rho \right) ,\varphi \in \left(
0,2\pi \right) ,\theta \in \left( 0,\pi \right) , 
\]%
\[
x_{1}=r\cos \varphi \sin \theta ,x_{2}=r\sin \varphi \sin \theta
,x_{3}=r\cos \theta . 
\]%
Also, denote%
\[
S_{\rho }=\left\{ r=\rho \right\} ,S_{\mu }=\left\{ r=\mu \right\} . 
\]%
The Laplace operator in the spherical coordinates is%
\begin{equation}
\Delta _{\text{sp}}w=w_{rr}+\frac{1}{r^{2}\sin ^{2}\theta }w_{\varphi
\varphi }+\frac{1}{\sin \theta }\frac{\partial }{\partial \theta }\left(
\sin \theta w_{\theta }\right) +\frac{2}{r}w_{r}=\widehat{\Delta }_{\text{sp}%
}w+\frac{2}{r}w_{r},  \label{3.10}
\end{equation}%
\begin{equation}
\widehat{\Delta }_{\text{sp}}w=w_{rr}+\frac{1}{r^{2}\sin ^{2}\theta }%
w_{\varphi \varphi }+\frac{1}{r^{2}\sin \theta }\frac{\partial }{\partial
\theta }\left( \sin \theta w_{\theta }\right) ,  \label{3.11}
\end{equation}%
for an arbitrary function $w\in C^{2}\left( \overline{\Omega }_{\mu }\right) 
$. We single out the operator $\widehat{\Delta }_{\text{sp}}$ in (\ref{3.10}%
), (\ref{3.11}) since any Carleman estimate is independent on the low order
derivatives of an operator, and also since we work in $\Omega _{\mu }$ where 
$r>\mu >0.$ Everywhere below $C=C\left( \Omega _{\mu }\right) >0$ denotes
different constants depending only on the domain $\Omega _{\mu }.$ Let 
\[
\nabla w=\left( w_{x_{1}},w_{x_{2}},w_{x_{3}}\right) ^{T}\text{ and }\nabla
_{\text{sp}}w=\left( w_{r},\frac{w_{\varphi }}{r\sin \theta },\frac{%
w_{\theta }}{r}\right) ^{T}. 
\]%
Note that since $w_{\varphi }=-w_{x_{1}}\sin \varphi \sin \theta
+w_{x_{2}}\cos \varphi \sin \theta ,$ then the function $w_{\varphi }/\sin
\theta $ does not have a singularity. It is well known that 
\begin{equation}
\left\vert \nabla w\right\vert \leq C\left\vert \nabla _{\text{sp}%
}w\right\vert \text{ in }\overline{\Omega }_{\mu },  \label{3.110}
\end{equation}%
\begin{equation}
\left\vert \nabla _{\text{sp}}w\right\vert \leq C\left\vert \nabla
w\right\vert \text{ in }\overline{\Omega }_{\mu }.  \label{3.111}
\end{equation}%
Introduce the subspace $H_{0}^{m}\left( \Omega _{\mu }\right) $ of the
Hilbert space $H^{m}\left( \Omega _{\mu }\right) $ as%
\[
H_{0}^{m}\left( \Omega _{\mu }\right) =\left\{ u\in H^{m}\left( \Omega _{\mu
}\right) :u\mid _{S_{\rho }}=u_{r}\mid _{S_{\rho }}=0\right\} ,\text{ }%
m=2,3. 
\]

We include the term with $\left( \Delta w\right) ^{2}$ in the right hand
side of the Carleman estimate (\ref{3.12}) since we will need to estimate
not only convergence for the vector function $W\left( x\right) $ (Theorem
5.4), but also to estimate convergence for the target coefficient $a_{0}(x)$
(Theorem 5.5). To do the latter, we will need to use equation (\ref{eqn:eitv}%
) in which the Laplace operator is involved.

\textbf{Theorem 3.1} (Carleman estimate). \emph{There exists a number }$%
\lambda _{0}=\lambda _{0}\left( \Omega _{\mu }\right) \geq 1$\emph{\ and a
number }$C=C\left( \Omega _{\mu }\right) >0,$\emph{\ both depending only on
the domain }$\Omega _{\mu },$\emph{\ such that for any function }$w\in
H^{2}\left( \Omega _{\mu }\right) $\emph{\ and for all }$\lambda \geq
\lambda _{0}$\emph{\ the following Carleman estimate with the CWF }$%
e^{2\lambda r}$\emph{\ holds:}%
\begin{equation}
\mathop{\displaystyle \int}\limits_{\Omega _{\mu }}\left( \Delta w\right)
^{2}e^{2\lambda r}dx\geq \frac{1}{2}\mathop{\displaystyle \int}%
\limits_{\Omega _{\mu }}\left( \Delta w\right) ^{2}e^{2\lambda r}dx+C\lambda %
\mathop{\displaystyle \int}\limits_{\Omega _{\mu }}\left( \nabla w\right)
^{2}e^{2\lambda r}dx+C\lambda ^{3}\mathop{\displaystyle \int}\limits_{\Omega
_{\mu }}w^{2}e^{2\lambda r}dx  \label{3.12}
\end{equation}%
\[
-C\lambda e^{2\lambda \rho }\mathop{\displaystyle \int}\limits_{S_{\rho
}}w_{r}^{2}dS-C\lambda ^{3}e^{2\lambda \rho }\mathop{\displaystyle \int}%
\limits_{S_{\rho }}w^{2}dS-C\lambda ^{3}e^{2\lambda \mu }\left\Vert
w\right\Vert _{H^{2}\left( \Omega _{\mu }\right) }^{2}. 
\]%
\emph{In particular, if }$w\in H_{0}^{2}\left( \Omega _{\mu }\right) ,$\emph{%
\ then}%
\begin{equation}
\mathop{\displaystyle \int}\limits_{\Omega _{\mu }}\left( \Delta w\right)
^{2}e^{2\lambda r}dx\geq \frac{1}{2}\mathop{\displaystyle \int}%
\limits_{\Omega _{\mu }}\left( \Delta w\right) ^{2}e^{2\lambda r}dx+C\lambda %
\mathop{\displaystyle \int}\limits_{\Omega _{\mu }}\left( \nabla w\right)
^{2}e^{2\lambda r}dx+C\lambda ^{3}\mathop{\displaystyle \int}\limits_{\Omega
_{\mu }}w^{2}e^{2\lambda r}dx  \label{3.122}
\end{equation}%
\[
-C\lambda ^{3}e^{2\lambda \mu }\left\Vert w\right\Vert _{H^{2}\left( \Omega
_{\mu }\right) }^{2}. 
\]

\textbf{Proof}. We assume that $w\in C^{2}\left( \overline{\Omega }_{\mu
}\right) $ since the case $w\in H^{2}\left( \Omega _{\mu }\right) $ can be
handled automatically later via density arguments. Introduce the new
function $q=we^{\lambda r}.$ Then%
\[
w=qe^{-\lambda r},w_{rr}=\left( q_{rr}-2\lambda q_{r}+\lambda ^{2}q\right)
e^{-\lambda r}. 
\]%
By (\ref{3.11})%
\[
\left( \widehat{\Delta }_{\text{sp}}w\right) ^{2}e^{2\lambda r}r\sin \theta =%
\left[ \left( \widehat{\Delta }_{\text{sp}}q+\lambda ^{2}q\right) -2\lambda
q_{r}\right] ^{2}r\sin \theta 
\]%
\[
\geq -4\lambda q_{r}\left( r\sin \theta q_{rr}+\frac{1}{r\sin \theta }%
q_{\varphi \varphi }+\frac{1}{r}\frac{\partial }{\partial \theta }\left(
\sin \theta q_{\theta }\right) +\lambda ^{2}r\sin \theta q\right) 
\]%
\[
=\partial _{r}\left( -2\lambda r\sin \theta q_{r}^{2}\right) +2\lambda \sin
\theta q_{r}^{2}+\partial _{\varphi }\left( -4\lambda \frac{q_{r}q_{\varphi }%
}{r\sin \theta }\right) +4\lambda \frac{1}{r\sin \theta }q_{r\varphi
}q_{\varphi } 
\]%
\[
+\partial _{\theta }\left( -4\lambda \sin \theta \frac{q_{r}q_{\theta }}{r}%
\right) +4\lambda \sin \theta \frac{q_{r\theta }q_{\theta }}{r}+\partial
_{r}\left( -2\lambda ^{3}r\sin \theta q^{2}\right) +2\lambda ^{3}r\sin
\theta q^{2} 
\]%
\[
=\partial _{r}\left( -2\lambda r\sin \theta q_{r}^{2}-2\lambda ^{3}r\sin
\theta q^{2}+2\lambda \frac{1}{r\sin \theta }q_{\varphi }^{2}+2\lambda \frac{%
\sin \theta }{r}q_{\theta }^{2}\right) 
\]%
\[
+\partial _{\varphi }\left( -4\lambda \frac{q_{r}q_{\varphi }}{r\sin \theta }%
\right) +\partial _{\theta }\left( -4\lambda \sin \theta \frac{%
q_{r}q_{\theta }}{r}\right) 
\]%
\[
+2\lambda \left( \sin \theta q_{r}^{2}+\frac{q_{\varphi }^{2}}{\sin \theta
r^{2}}+\sin \theta \frac{q_{\theta }^{2}}{r^{2}}\right) +\lambda ^{3}r\sin
\theta q^{2}. 
\]%
Hence, we have proven that%
\[
\left( \widehat{\Delta }_{\text{sp}}w\right) ^{2}e^{2\lambda r}r\sin \theta 
\]%
\[
\geq 2\lambda \left( \sin \theta q_{r}^{2}+\frac{q_{\varphi }^{2}}{\sin
\theta r^{2}}+\sin \theta \frac{q_{\theta }^{2}}{r^{2}}\right) +2\lambda
^{3}r\sin \theta q^{2} 
\]%
\[
+\partial _{r}\left( -2\lambda r\sin \theta q_{r}^{2}-2\lambda ^{3}r\sin
\theta q^{2}+2\lambda \frac{1}{r\sin \theta }q_{\varphi }^{2}+2\lambda \frac{%
\sin \theta }{r}q_{\theta }^{2}\right) 
\]%
\[
+\partial _{\varphi }\left( -4\lambda \frac{q_{r}q_{\varphi }}{r\sin \theta }%
\right) +\partial _{\theta }\left( -4\lambda \sin \theta \frac{%
q_{r}q_{\theta }}{r}\right) . 
\]%
Integrate this inequality over $\Omega _{\mu }$ while keeping in mind that
the function $q\left( r,\varphi ,\theta \right) $ is periodic with respect
to $\varphi $ with the period $2\pi $ and that $\sin 0=\sin \pi =0$ and also
that $dx=r^{2}\sin \theta drd\varphi d\theta .$ We obtain%
\[
\mathop{\displaystyle \int}\limits_{\Omega _{\mu }}\left( \widehat{\Delta }_{%
\text{sp}}w\right) ^{2}e^{2\lambda r}dx\geq C\mathop{\displaystyle \int}%
\limits_{\Omega _{\mu }}\left( \widehat{\Delta }_{\text{sp}}w\right) ^{2}%
\frac{e^{2\lambda r}}{r}dx=\mathop{\displaystyle \int}\limits_{\Omega _{\mu
}}\left( \widehat{\Delta }_{\text{sp}}w\right) ^{2}e^{2\lambda r}r\sin
\theta drd\varphi d\theta 
\]%
\begin{equation}
\geq 2\lambda \mathop{\displaystyle \int}\limits_{\Omega _{\mu }}\left(
q_{r}^{2}+\frac{q_{\varphi }^{2}}{\sin ^{2}\theta r^{2}}+\frac{q_{\theta
}^{2}}{r^{2}}\right) \sin \theta drd\varphi d\theta +2\lambda ^{3}%
\mathop{\displaystyle \int}\limits_{\Omega _{\mu }}q^{2}r\sin \theta
drd\varphi d\theta  \label{3.121}
\end{equation}%
\[
-\mathop{\displaystyle \int}\limits_{S_{\rho }}\left( 2\lambda
q_{r}^{2}+2\lambda ^{3}q^{2}\right) dS-C\lambda e^{2\lambda \mu }%
\mathop{\displaystyle \int}\limits_{S_{\mu }}\left( \nabla w\right) ^{2}dS. 
\]%
Change variables back from $q$ to $w$. Since $q=we^{\lambda r},$ then%
\[
q_{r}^{2}=\left( w_{r}+\lambda w\right) ^{2}e^{2\lambda
r}=w_{r}^{2}e^{2\lambda r}+2\lambda w_{r}we^{2\lambda r}+\lambda
^{2}w^{2}e^{2\lambda r} 
\]%
\[
=w_{r}^{2}e^{2\lambda r}+\partial _{r}\left( \lambda w^{2}e^{2\lambda
r}\right) -2\lambda ^{2}w^{2}e^{2\lambda r}+\lambda ^{2}w^{2}e^{2\lambda r} 
\]%
\[
=w_{r}^{2}e^{2\lambda r}-\lambda ^{2}w^{2}e^{2\lambda r}+\partial _{r}\left(
\lambda w^{2}e^{2\lambda r}\right) . 
\]%
Let the number $a=\min \left( \mu /2,1\right) .$ Then in the first line of (%
\ref{3.121})%
\[
2\lambda q_{r}^{2}+2\lambda ^{3}q^{2}r\geq 2\lambda aq_{r}^{2}+2\lambda
^{3}q^{2}\mu 
\]%
\begin{equation}
\geq 2\lambda aw_{r}^{2}e^{2\lambda r}+2\lambda ^{3}\left( \mu -\frac{\mu }{2%
}\right) w^{2}e^{2\lambda r}+\partial _{r}\left( 2\lambda
^{2}aw^{2}e^{2\lambda r}\right)  \label{3.123}
\end{equation}%
\[
\geq C\lambda w_{r}^{2}e^{2\lambda r}+C\lambda ^{3}w^{2}e^{2\lambda
r}+\partial _{r}\left( 2\lambda ^{2}aw^{2}e^{2\lambda r}\right) . 
\]%
Hence, using (\ref{3.110}), (\ref{3.121}) and (\ref{3.123}), we obtain%
\begin{equation}
\mathop{\displaystyle \int}\limits_{\Omega _{\mu }}\left( \widehat{\Delta }_{%
\text{sp}}w\right) ^{2}e^{2\lambda r}dx\geq C\lambda \mathop{\displaystyle
\int}\limits_{\Omega _{\mu }}\left( \nabla w\right) ^{2}e^{2\lambda
r}dx+C\lambda ^{3}\mathop{\displaystyle \int}\limits_{\Omega _{\mu
}}w^{2}e^{2\lambda r}dx  \label{3.124}
\end{equation}%
\[
-C\lambda e^{2\lambda \rho }\mathop{\displaystyle \int}\limits_{S_{\rho
}}w_{r}^{2}dS-C\lambda ^{3}e^{2\lambda \rho }\mathop{\displaystyle \int}%
\limits_{S_{\rho }}w^{2}dS-C\lambda ^{3}e^{2\lambda \mu }\mathop{%
\displaystyle \int}\limits_{S_{\mu }}\left( \left( \nabla w\right)
^{2}+w^{2}\right) dS. 
\]%
Noticing that by (\ref{3.10}) 
\[
\left( \Delta w\right) ^{2}=\left( \widehat{\Delta }_{\text{sp}}w\right)
^{2}+4\left( \widehat{\Delta }_{\text{sp}}w\right) \frac{w_{r}}{r}+4\left( 
\frac{w_{r}}{r}\right) ^{2}\geq \frac{1}{2}\left( \widehat{\Delta }_{\text{sp%
}}w\right) ^{2}-Cw_{r}^{2}, 
\]%
and also that 
\[
\mathop{\displaystyle \int}\limits_{S_{\mu }}\left( \left( \nabla w\right)
^{2}+w^{2}\right) dS\leq C\left\Vert w\right\Vert _{H^{2}\left( \Omega _{\mu
}\right) }^{2}, 
\]%
and then using (\ref{3.124}), we obtain 
\begin{equation}
\mathop{\displaystyle \int}\limits_{\Omega _{\mu }}\left( \Delta w\right)
^{2}e^{2\lambda r}dx\geq C\lambda \mathop{\displaystyle \int}\limits_{\Omega
_{\mu }}\left( \nabla w\right) ^{2}e^{2\lambda r}dx+C\lambda ^{3}%
\mathop{\displaystyle \int}\limits_{\Omega _{\mu }}w^{2}e^{2\lambda r}dx
\label{3.125}
\end{equation}%
\[
-C\lambda e^{2\lambda \rho }\mathop{\displaystyle \int}\limits_{S_{\rho
}}w_{r}^{2}dS-C\lambda ^{3}e^{2\lambda \rho }\mathop{\displaystyle \int}%
\limits_{S_{\rho }}w^{2}dS-C\lambda ^{3}e^{2\lambda \mu }\left\Vert
w\right\Vert _{H^{2}\left( \Omega _{\mu }\right) }^{2}. 
\]%
Obviously,%
\begin{equation}
\mathop{\displaystyle \int}\limits_{\Omega _{\mu }}\left( \Delta w\right)
^{2}e^{2\lambda r}dx=\mathop{\displaystyle \int}\limits_{\Omega _{\mu
}}\left( \Delta w\right) ^{2}e^{2\lambda r}dx.  \label{3.126}
\end{equation}%
Summing up (\ref{3.125}) with (\ref{3.126}) and then dividing the resulting
estimate by 2, we obtain (\ref{3.12}). $\square $

\subsubsection{The 2D case}

\label{sect:3.3.2}

In this case we keep the same notations for domains $\Omega ,\Omega _{\mu }$
as ones in subsection 3.3.1, meaning, however, that now these are domains in 
$\mathbb{R}^{2}.$ Polar coordinates are 
\[
r=\left\vert x\right\vert \in \left( \mu ,\rho \right) ,\varphi \in \left(
0,2\pi \right) , 
\]%
\[
x_{1}=r\cos \varphi ,x_{2}=r\sin \varphi . 
\]%
The Laplace operator in polar coordinates is%
\[
\Delta _{p}w=w_{rr}+\frac{1}{r^{2}}w_{\varphi \varphi }+\frac{1}{r}w_{r}=%
\widehat{\Delta }_{\text{p}}w+\frac{1}{r}w_{r}. 
\]%
\textbf{Theorem 3.2} (Carleman estimate). \emph{There exists a number }$%
\lambda _{0}=\lambda _{0}\left( \Omega _{\mu }\right) \geq 1$\emph{\
depending only on the domain }$\Omega _{\mu }$\emph{\ such that for any
function }$w\in H^{2}\left( \Omega _{\mu }\right) $\emph{\ and for all }$%
\lambda \geq \lambda _{0}$\emph{\ the following Carleman estimate holds:}%
\[
\mathop{\displaystyle \int}\limits_{\Omega _{\mu }}\left( \Delta w\right)
^{2}e^{2\lambda r}dx\geq \frac{1}{2}\mathop{\displaystyle \int}%
\limits_{\Omega _{\mu }}\left( \Delta w\right) ^{2}e^{2\lambda r}dx+C\lambda %
\mathop{\displaystyle \int}\limits_{\Omega _{\mu }}\left( \nabla w\right)
^{2}e^{2\lambda r}dx+C\lambda ^{3}\mathop{\displaystyle \int}\limits_{\Omega
_{\mu }}w^{2}e^{2\lambda r}dx 
\]%
\[
-C\lambda e^{2\lambda \rho }\mathop{\displaystyle \int}\limits_{S_{\rho
}}w_{r}^{2}dS-C\lambda ^{3}e^{2\lambda \rho }\mathop{\displaystyle \int}%
\limits_{S_{\rho }}w^{2}dS-C\lambda e^{2\lambda \mu }\left\Vert w\right\Vert
_{H^{2}\left( \Omega _{\mu }\right) }^{2}. 
\]%
In particular, if $w\in H_{0}^{2}\left( \Omega _{\mu }\right) ,$\emph{\ then}%
\[
\mathop{\displaystyle \int}\limits_{\Omega _{\mu }}\left( \Delta w\right)
^{2}e^{2\lambda r}dx\geq \frac{1}{2}\mathop{\displaystyle \int}%
\limits_{\Omega _{\mu }}\left( \Delta w\right) ^{2}e^{2\lambda r}dx+C\lambda %
\mathop{\displaystyle \int}\limits_{\Omega _{\mu }}\left( \nabla w\right)
^{2}e^{2\lambda r}dx+C\lambda ^{3}\mathop{\displaystyle \int}\limits_{\Omega
_{\mu }}w^{2}e^{2\lambda r}dx 
\]%
\[
-C\lambda ^{3}e^{2\lambda \mu }\left\Vert w\right\Vert _{H^{2}\left( \Omega
_{\mu }\right) }^{2}. 
\]

The proof of this theorem is omitted since it is very similar with the proof
of Theorem 3.1.

\subsection{H\"{o}lder stability and uniqueness of the Cauchy problem (%
\protect\ref{3.8}), (\protect\ref{3.9})}

\label{sect:3.4}

We establish in this subsection the H\"{o}lder stability estimate for
problem (\ref{3.8}), (\ref{3.9}). Uniqueness follows immediately from this
estimate. We work here only with the 3D case. Theorem 3.2 implies that the
2D case can be handled almost exactly the same way. Thus, in this subsection
the domain $\Omega $ is as in (\ref{3.900}), and in (\ref{3.9}) $\Gamma
=\partial \Omega =\left\{ r=\rho \right\} .$ Everywhere below we often work
with $N-$dimensional vector functions, like, e.g. $V\left( x\right) .$ Norms
in standard functional spaces of such vector functions are defined in the
natural well known way via corresponding norms of their components. The same
about scalar products. It is always clear from the context whether we work
with regular functions or with those $N-$dimensional vector functions.

Suppose that there exist two solutions of problem (\ref{3.8}), (\ref{3.9}), $%
V_{1},V_{2}\in H^{2}\left( \Omega \right) \cap C^{1}\left( \overline{\Omega }%
\right) $ such that 
\begin{equation}
V_{1}\mid _{S_{\rho }}=p_{0},V_{2}\mid _{S_{\rho }}=p_{0,\delta
},V_{1,r}\mid _{S_{\rho }}=p_{1},V_{2,r}\mid _{S_{\rho }}=p_{1,\delta },
\label{4.100}
\end{equation}%
where 
\begin{equation}
\left\Vert p_{0,\delta }-p_{0}\right\Vert _{L^{2}\left( S_{\rho }\right)
}\leq \delta ,\left\Vert p_{1,\delta }-p_{1}\right\Vert _{L^{2}\left(
S_{\rho }\right) }\leq \delta ,  \label{3.120}
\end{equation}%
where $\delta \in \left( 0,1\right) $ is a sufficiently small number which
is interpreted as the level of the noise in the data. Denote 
\begin{equation}
\widetilde{V}=V_{1}-V_{2},\widetilde{p}=p_{0}-p_{0,\delta },\widetilde{p}%
_{1}=p_{1}-p_{1,\delta }.  \label{4.101}
\end{equation}%
Recalling that the function $F$ in (\ref{3.8}) is quadratic with respect to
the derivatives $\partial _{x_{j}}v_{k}\left( x\right) ,$ we obtain from (%
\ref{3.8}) and (\ref{3.9})%
\begin{equation}
\Delta \widetilde{V}=Q\left( \nabla V_{1},\nabla V_{2}\right) \nabla 
\widetilde{V},\text{ }x\in \Omega ,\widetilde{V}\in H^{2}\left( \Omega
\right) \cap C^{1}\left( \overline{\Omega }\right) ,  \label{3.13}
\end{equation}%
\begin{equation}
\widetilde{V}\mid _{S_{\rho }}=\widetilde{p}_{0},\widetilde{V}_{r}\mid
_{S_{\rho }}=\widetilde{p}_{1},  \label{3.14}
\end{equation}%
where the vector function $Q\left( \nabla V_{1},\nabla V_{2}\right) $ is
linear with respect to components of vector functions $\nabla V_{1},\nabla
V_{2}.$

\textbf{Theorem 3.3} (H\"{o}lder stability estimate). \emph{For two vector
functions }$V_{1},V_{2}\in H^{2}\left( \Omega \right) \cap C^{1}\left( 
\overline{\Omega }\right) $\emph{\ introduced above in this section, let }$%
\left\Vert V_{1}\right\Vert _{C^{1}\left( \overline{\Omega }\right)
},\left\Vert V_{2}\right\Vert _{C^{1}\left( \overline{\Omega }\right) }\leq
A,$\emph{\ where }$A=const.>0.$\emph{\ Let (\ref{4.100})-(\ref{4.101}) hold.
Choose a number }$\eta \in \left( 0,\rho -\mu \right) .$\emph{\ Let }$\Omega
_{\mu +\eta }=\left\{ x:\mu +\eta <\left\vert x\right\vert <\rho \right\}
\subset \Omega _{\mu }.$\emph{\ Then there exists a number }$%
C_{1}=C_{1}\left( \Omega _{\mu },\eta ,F,\Psi \left( N\right) ,A\right) >0$%
\emph{\ and a sufficiently small number }$\delta _{0}=\delta _{0}\left(
\Omega _{\mu },\eta ,F,\Psi \left( N\right) ,A\right) \in \left( 0,1\right) $%
\emph{\ such that for all }$\delta \in \left( 0,\delta _{0}\right) $\emph{\
the following H\"{o}lder stability estimate holds:}%
\begin{equation}
\left\Vert \widetilde{V}\right\Vert _{H^{1}\left( \Omega _{\mu +\eta
}\right) }\leq C_{1}\delta ^{\gamma },\text{ where }\gamma =\eta /\left(
4\rho \right) .  \label{3.150}
\end{equation}

\textbf{Proof}. In this proof $C_{1}=C_{1}\left( \Omega _{\mu },\eta ,F,\Psi
\left( N\right) ,A\right) >0$ denotes different positive constants depending
only on listed parameters. Note that $\left\vert Q\left( \nabla V_{1},\nabla
V_{2}\right) \right\vert \leq C_{1}.$ A careful analysis of the proof of
Theorem 3.1, more precisely of the last term in the third line of (\ref%
{3.124}), shows that the term $\left\Vert w\right\Vert _{H^{2}\left( \Omega
_{\mu }\right) }^{2}$ in (\ref{3.12}) can be replaced with the term $%
\left\Vert w\right\Vert _{C^{1}\left( \overline{\Omega }_{\mu }\right)
}^{2}. $ Squaring both sides of (\ref{3.13}), replacing the equality with
the inequality and using (\ref{3.13}), we obtain%
\begin{equation}
\left( \Delta \widetilde{V}\right) ^{2}\leq C_{1}\left( \nabla \widetilde{V}%
\right) ^{2},\text{ }x\in \Omega _{\mu }.  \label{3.16}
\end{equation}%
Multiplying both sides of (\ref{3.16}) by $e^{2\lambda r}$ and integrating
over the domain $\Omega _{\mu },$ we obtain%
\begin{equation}
C_{1}\mathop{\displaystyle \int}\limits_{\Omega _{\mu }}\left( \nabla 
\widetilde{V}\right) ^{2}e^{2\lambda r}dx\geq \mathop{\displaystyle \int}%
\limits_{\Omega _{\mu }}\left( \Delta \widetilde{V}\right) ^{2}e^{2\lambda
r}dx.  \label{3.17}
\end{equation}%
Next, by (\ref{3.12}) )%
\[
\mathop{\displaystyle \int}\limits_{\Omega _{\mu }}\left( \Delta \widetilde{V%
}\right) ^{2}e^{2\lambda r}dx\geq C\lambda \mathop{\displaystyle \int}%
\limits_{\Omega _{\mu }}\left( \nabla \widetilde{V}\right) ^{2}e^{2\lambda
r}dx+C\lambda ^{3}\mathop{\displaystyle \int}\limits_{\Omega _{\mu }}%
\widetilde{V}^{2}e^{2\lambda r}dx 
\]%
\[
-C\lambda e^{2\lambda \rho }\mathop{\displaystyle \int}\limits_{S_{\rho }}%
\widetilde{p}_{1}^{2}dS-C\lambda ^{3}e^{2\lambda \rho }%
\mathop{\displaystyle
\int}\limits_{S_{\rho }}\widetilde{p}_{0}^{2}dS-C\lambda ^{3}e^{2\lambda \mu
}\left\Vert \widetilde{V}\right\Vert _{C^{1}\left( \overline{\Omega }_{\mu
}\right) }^{2} 
\]%
\[
\geq C\lambda \mathop{\displaystyle \int}\limits_{\Omega _{\mu }}\left(
\nabla \widetilde{V}\right) ^{2}e^{2\lambda r}dx+C\lambda ^{3}%
\mathop{\displaystyle \int}\limits_{\Omega _{\mu }}\widetilde{V}%
^{2}e^{2\lambda r}dx 
\]%
\[
-Ce^{2\lambda \rho }\delta ^{2}-Ce^{2\lambda \mu }\left\Vert \widetilde{V}%
\right\Vert _{C^{1}\left( \overline{\Omega }_{\mu }\right) }^{2}. 
\]%
Hence, taking into account (\ref{3.17}), we obtain for sufficiently large $%
\lambda _{1}=\lambda _{1}\left( C,C_{1}\right) \geq \lambda _{0}>0$%
\begin{equation}
C_{1}e^{2\lambda \rho }\delta ^{2}+C_{1}\lambda ^{3}e^{2\lambda \mu
}\left\Vert \widetilde{V}\right\Vert _{C^{1}\left( \overline{\Omega }_{\mu
}\right) }^{2}\geq \lambda \mathop{\displaystyle \int}\limits_{\Omega _{\mu
}}\left( \nabla \widetilde{V}\right) ^{2}e^{2\lambda r}dx+\lambda ^{3}%
\mathop{\displaystyle \int}\limits_{\Omega _{\mu }}\widetilde{V}%
^{2}e^{2\lambda r}dx.  \label{3.18}
\end{equation}%
Since $\Omega _{\mu +\eta }\subset \Omega _{\mu }$ and $e^{2\lambda
r}>e^{2\lambda \left( \mu +\eta \right) }$ in $\Omega _{\mu +\eta },$ then%
\[
\lambda \mathop{\displaystyle \int}\limits_{\Omega _{\mu }}\left( \nabla 
\widetilde{V}\right) ^{2}e^{2\lambda r}dx+\lambda ^{3}%
\mathop{\displaystyle
\int}\limits_{\Omega _{\mu }}\widetilde{V}^{2}e^{2\lambda r}dx\geq
e^{2\lambda \left( \mu +\eta \right) }\mathop{\displaystyle \int}%
\limits_{\Omega _{\mu +\eta }}\left( \left( \nabla \widetilde{V}\right) ^{2}+%
\widetilde{V}^{2}\right) dx. 
\]%
Comparing this with (\ref{3.18}), we obtain%
\begin{equation}
\left\Vert \widetilde{V}\right\Vert _{H^{1}\left( \Omega _{\mu +\eta
}\right) }^{2}\leq C_{1}e^{2\lambda \rho }\delta ^{2}+C_{1}\lambda
^{3}e^{-2\lambda \eta }\left\Vert \widetilde{V}\right\Vert _{C^{1}\left( 
\overline{\Omega }_{\mu }\right) }^{2}.  \label{3.19}
\end{equation}%
Since $\lambda ^{3}e^{-2\lambda \eta }\leq e^{-\lambda \eta }$ for
sufficiently large $\lambda \geq \lambda \left( C,C_{1}\right) >1,$ then (%
\ref{3.19}) implies that 
\begin{equation}
\left\Vert \widetilde{V}\right\Vert _{H^{1}\left( \Omega _{\mu +\eta
}\right) }^{2}\leq C_{1}e^{2\lambda \rho }\delta ^{2}+C_{1}e^{-\lambda \eta
}\left\Vert \widetilde{V}\right\Vert _{C^{1}\left( \overline{\Omega }_{\mu
}\right) }^{2}.  \label{3.20}
\end{equation}%
Choose $\lambda =\lambda \left( \delta \right) $ such that $e^{2\lambda \rho
}\delta ^{2}=\delta .$ Hence, $\lambda =\ln \left( \delta ^{-1/\left( 2\rho
\right) }\right) .$ Since we must have $\lambda \geq \lambda _{1}=\lambda
_{1}\left( C,C_{1}\right) \geq \lambda _{0}>0,$ then we must have $\delta
<\delta _{0}=\exp \left( -2\rho \lambda _{1}\right) .$ Next, $e^{-\lambda
\eta }=\delta ^{\eta /\left( 2\rho \right) }.$ Since $\eta /\left( 2\rho
\right) <1/2,$ then $\delta ^{\eta /\left( 2\rho \right) }>\delta .$ Set $%
2\kappa =\eta /\left( 2\rho \right) \in \left( 0,1/2\right) .$ Noticing that 
$\left\Vert \widetilde{V}\right\Vert _{C^{1}\left( \overline{\Omega }_{\mu
}\right) }^{2}\leq 2\left( \left\Vert V_{1}\right\Vert _{C^{1}\left( 
\overline{\Omega }_{\mu }\right) }^{2}+\left\Vert V_{2}\right\Vert
_{C^{1}\left( \overline{\Omega }_{\mu }\right) }^{2}\right) \leq 2A^{2},$ we
obtain from (\ref{3.20}) the target estimate (\ref{3.150}). \ $\square $

\section{Convexification}

\label{sec:4}

To solve the Cauchy problem (\ref{3.8}), (\ref{3.9}) numerically, we
construct in this section a weighted Tikhonov-like functional with the CWF $%
e^{2\lambda r}$ in it and prove necessary theorems. For brevity, we
construct the Tikhonov-like functional only for the 3D case. So, in sections
4 and 5 we work only with the 3D case. The 2D case is completely similar and
direct analogs of Theorems 5.1-5.4 (below) are valid in 2D.

\subsection{Weighted Tikhonov-like functional}

\label{sec:4.1}

We assume that in (\ref{3.9})%
\begin{equation}
\Gamma =\partial \Omega =S_{\rho };\text{ }p_{0},p_{1}\in C^{3}\left(
S_{\rho }\right) .  \label{4.1}
\end{equation}%
We now arrange zero Dirichlet and Neumann boundary conditions for a new
vector function $W$, which is associated with the vector function $V$. We
are doing so since we use below some theorems of \cite{BakKlib}, which are
applicable only in the case of zero boundary conditions.

Denote%
\begin{equation}
P\left( r,\varphi ,\theta \right) =p_{0}\left( r,\varphi ,\theta \right)
+\left( r-\rho \right) p_{1}\left( r,\varphi ,\theta \right) ,  \label{4.3}
\end{equation}%
\begin{equation}
W\left( r,\varphi ,\theta \right) =V\left( r,\varphi ,\theta \right)
-P\left( r,\varphi ,\theta \right) ;\text{ }W\left( r,\varphi ,\theta
\right) =\left( W_{0},...,W_{N-1}\right) ^{T}\left( r,\varphi ,\theta
\right) .  \label{4.4}
\end{equation}%
Then by (\ref{4.1}) $P\in C^{3}\left( \overline{\Omega }_{\mu }\right) $.
Hence, (\ref{3.8}), (\ref{3.9}), (\ref{4.3}) and (\ref{4.4}) imply that%
\begin{equation}
\Delta W+\Delta P-F\left( \nabla W+\nabla P\right) =0,\text{ }  \label{4.5}
\end{equation}%
\begin{equation}
W\in H_{0}^{3}\left( \Omega _{\mu }\right) .  \label{4.6}
\end{equation}%
Note that by the embedding theorem 
\begin{equation}
H^{3}\left( \Omega _{\mu }\right) \subset C^{1}\left( \overline{\Omega }%
_{\mu }\right) ,\left\Vert f\right\Vert _{C^{1}\left( \overline{\Omega }%
_{\mu }\right) }\leq C\left\Vert f\right\Vert _{H^{3}\left( \Omega _{\mu
}\right) }.  \label{4.8}
\end{equation}

Let $\eta \in \left( 0,\rho -\mu \right) $ be the number which was chosen in
Theorem 3.3. Our weighted Tikhonov-like functional is: 
\begin{equation}
J_{\lambda ,\beta }\left( W\right) =  \label{4.9}
\end{equation}%
\[
=e^{-2\lambda \left( \mu +\eta \right) }\mathop{\displaystyle \int}%
\limits_{\Omega _{\mu }}\left[ \Delta W+\Delta P-F\left( \nabla W+\nabla
P\right) \right] ^{2}e^{2\lambda r}dx+\beta \left\Vert W+P\right\Vert
_{H^{3}\left( \Omega _{\mu }\right) }^{2}. 
\]%
Here $\beta \in \left( 0,1\right) $ is the regularization parameter and the
multiplier $e^{-2\lambda \left( \mu +\eta \right) }$ is introduced to
balance two terms in the right hand side of (\ref{4.9}). Let $R>0$ be an
arbitrary number. Let $B\left( R\right) \subset H_{0}^{3}\left( \Omega _{\mu
}\right) $ be the ball of the radius $R$ with the center at $\left\{
0\right\} ,$ 
\begin{equation}
B\left( R\right) =\left\{ W\in H_{0}^{3}\left( \Omega _{\mu }\right)
:\left\Vert W\right\Vert _{H^{3}\left( \Omega _{\mu }\right) }<R\right\} .
\label{4.10}
\end{equation}%
We consider the following minimization problem:

\textbf{Minimization Problem}. \emph{Minimize the functional }$J_{\lambda
,\beta }\left( W\right) $\emph{\ on the closed ball} $\overline{B\left(
R\right) }.$

\section{Theorems}

\label{sec:5}

In this section we formulate and prove some theorems about the above
minimization problem.

\subsection{Formulations of theorems}

\label{sec:5.1}

The central analytical result of this paper is Theorem 5.1.

\textbf{Theorem 5.1}. \emph{The functional }$J_{\lambda ,\beta }\left(
W\right) $\emph{\ has the Frech\'{e}t derivative }$J_{\lambda ,\beta
}^{\prime }\left( W\right) $\emph{\ at every point }$W\in H_{0}^{3}\left(
\Omega _{\mu }\right) .$\emph{\ Furthermore, there exists numbers }

$\lambda _{2}=\lambda _{2}\left( \mu ,\eta ,F,\Psi \left( N\right)
,P,R\right) \geq \lambda _{0}>0$\emph{\ and }$C_{2}=C_{2}\left( \mu ,\eta
,F,\Psi \left( N\right) ,P,R\right) >0$ \emph{depending only on listed
parameters such that }$2e^{-\lambda _{2}\eta }<1$\emph{\ and for all }$%
\lambda \geq \lambda _{2}$\emph{\ the functional }$J_{\lambda ,\beta }\left(
W\right) $\emph{\ is strictly convex on }$\overline{B\left( R\right) }$ 
\emph{for the choice of }$\beta $ \emph{as} \emph{\ }%
\begin{equation}
\beta \in \left( 2e^{-\lambda \eta },1\right) .  \label{5.0}
\end{equation}%
\emph{\ More precisely, the following inequality holds:}%
\[
J_{\lambda ,\beta }\left( W_{2}\right) -J_{\lambda ,\beta }\left(
W_{1}\right) -J_{\lambda ,\beta }^{\prime }\left( W_{1}\right) \left(
W_{2}-W_{1}\right) 
\]%
\begin{equation}
\geq C_{2}\left\Vert \Delta \left( W_{2}-W_{1}\right) \right\Vert
_{L^{2}\left( \Omega _{\mu +\eta }\right) }+C_{2}\left\Vert
W_{2}-W_{1}\right\Vert _{H^{1}\left( \Omega _{\mu +\eta }\right) }^{2}+\frac{%
\beta }{2}\left\Vert W_{2}-W_{1}\right\Vert _{H^{3}\left( \Omega _{\mu
}\right) }^{2},  \label{5.1}
\end{equation}%
\[
\text{ }\forall W_{1},W_{2}\in \overline{B\left( R\right) }. 
\]

Note that, allowing the regularization parameter $\beta \in \left(
2e^{-\lambda \eta },1\right) ,$ we actually allow $\beta $ to be
sufficiently small. We now formulate the theorem about the Lipschitz
continuity condition of the Frech\'{e}t derivative\emph{\ }$J_{\lambda
,\beta }^{\prime }\left( W\right) .$

\textbf{Theorem 5.2}. \emph{For any numbers }$\widetilde{R},\lambda >0,\beta
\in \left( 0,1\right) $\emph{\ the Frech\'{e}t derivative }$J_{\lambda
,\beta }^{\prime }\left( W\right) $\emph{\ of the functional }$J_{\lambda
,\beta }\left( W\right) $\emph{\ satisfies the Lipschitz continuity
condition in the ball }$B\left( \widetilde{R}\right) .$\emph{\ In other
words, there exists a number }$Z=Z\left( \Omega _{\mu },F,\Psi \left(
N\right) ,\widetilde{R},\lambda \right) >0$\emph{\ depending only on listed
parameters such that} 
\[
\left\vert J_{\lambda ,\beta }^{\prime }\left( W_{2}\right) -J_{\lambda
,\beta }^{\prime }\left( W_{1}\right) \right\vert \leq Z\left\Vert
W_{2}-W_{1}\right\Vert _{H^{3}\left( \Omega _{\mu }\right) },\text{ }\forall
W_{1},W_{2}\in B\left( \widetilde{R}\right) . 
\]

Consider now the gradient projection method of the minimization of the
functional $J_{\lambda ,\beta }$ on the closed ball $\overline{B\left(
R\right) }.$ Let $P_{B}:H_{0}^{3}\left( \Omega _{\mu }\right) \rightarrow 
\overline{B\left( R\right) }$ be the projection operator of the space $%
H_{0}^{3}\left( \Omega _{\mu }\right) $ on the closed ball $W_{\min ,\lambda
}\subset H_{0}^{3}\left( \Omega _{\mu }\right) .$ Let $W_{0}\in B\left(
R\right) $ be an arbitrary point. The sequence $\left\{ W_{n}\right\}
_{n=1}^{\infty }$ of the gradient projection method is defined as%
\begin{equation}
W_{n}=P_{B}\left( W_{n-1}-\zeta J_{\lambda ,\beta }^{\prime }\left(
W_{n-1}\right) \right) ,\text{ }n=1,2,...,  \label{5.2}
\end{equation}%
where $\zeta \in \left( 0,1\right) $ is a sufficiently small number. \ Below 
$\left[ ,\right] $ denotes the scalar product in the space of real valued $%
N- $D vector functions $H^{3}\left( \Omega _{\mu }\right) .$

\textbf{Theorem 5.3}. \emph{Let }$\lambda _{2}=\lambda _{2}\left( \mu ,\eta
,F,N,P,R\right) \geq \lambda _{0}>0$\emph{\ be the number of Theorem 5.1 and
let the regularization parameter }$\beta \in \left( 2e^{-\lambda \eta
},1\right) $\emph{. Then for every }$\lambda \geq \lambda _{2}$\emph{\ there
exists unique minimizer }$W_{\min ,\lambda ,\beta }\in \overline{B\left(
R\right) }$\emph{\ of the functional }$J_{\lambda ,\beta }\left( W\right) $%
\emph{\ on the closed ball }$\overline{B\left( R\right) }.$\emph{\
Furthermore, the following inequality holds}%
\begin{equation}
\left[ J_{\lambda ,\beta }^{\prime }\left( W_{\min ,\lambda ,\beta }\right)
,W-W_{\min ,\lambda ,\beta }\right] \geq 0,\text{ }\forall W\in \overline{%
B\left( R\right) }.  \label{5.200}
\end{equation}%
\emph{In addition, there exists a sufficiently small number }$\zeta
_{0}=\zeta _{0}\left( \mu ,\eta ,F,\Psi \left( N\right) ,P,R,\lambda ,\beta
\right) \in \left( 0,1\right) $\emph{\ depending only on listed parameters
such that for every }$\zeta \in \left( 0,\zeta _{0}\right) $\emph{\ the
sequence (\ref{5.2}) converges to the minimizer }$W_{\min ,\lambda ,\beta }$%
\emph{\ and the following estimate of the convergence rate holds:}%
\begin{equation}
\left\Vert W_{n}-W_{\min ,\lambda ,\beta }\right\Vert _{H^{3}\left( \Omega
\right) }\leq \omega ^{n}\left\Vert W_{0}-W_{\min ,\lambda ,\beta
}\right\Vert _{H^{3}\left( \Omega \right) },\text{ }n=1,2,...,  \label{5.201}
\end{equation}%
\emph{where the number }$\omega =\omega \left( \zeta \right) \in \left(
0,1\right) $\emph{\ depends only on the parameter }$\zeta .$

Even though Theorem 5.3 guarantees the convergence of the gradient
projection method to the unique minimizer of the functional (\ref{4.9}), it
is not yet clear how far this minimizer is from the exact solution. To
address this question, we assume, as it is commonly accepted in the theory
of ill-posed problems \cite{T}, that there exists an exact solution $W^{\ast
}\in B\left( R\right) $ of the problem (\ref{4.5}), (\ref{4.6}), i.e.
solution with the noiseless data.

Let $\delta \in \left( 0,1\right) $ be a sufficiently small number
characterizing the level of the noise in the data. Let $W^{\ast }$ be the
exact solution of problem (\ref{4.5}), (\ref{4.6}) with the noiseless data $%
P^{\ast }\in C^{3}\left( \overline{\Omega }_{\mu }\right) ,$%
\begin{equation}
\Delta W^{\ast }+\Delta P^{\ast }-F\left( \nabla W^{\ast }+\nabla P^{\ast
}\right) =0,\text{ }  \label{5.3}
\end{equation}%
\begin{equation}
W^{\ast }\in H_{0}^{3}\left( \Omega _{\mu }\right) .  \label{5.4}
\end{equation}%
Let $P\in C^{3}\left( \overline{\Omega }_{\mu }\right) $ be the noisy data.
Denote $\widetilde{P}=P-P^{\ast }.$ We assume that%
\begin{equation}
\left\Vert \widetilde{P}\right\Vert _{H^{3}\left( \overline{\Omega }_{\mu
}\right) }\leq \delta .  \label{5.5}
\end{equation}

\textbf{Theorem 5.4.} \emph{Let }$\lambda _{2}\geq \lambda _{0}>0$\emph{\
and }$C_{2}>0$\emph{\ be numbers of Theorem 5.1. Choose the number }$\delta
_{1}>0$\emph{\ so small that }$\delta _{1}<\min \left( e^{-4\rho \lambda
_{2}},3^{-4\rho /\eta }\right) $ \emph{and let} $\delta \in \left( 0,\delta
_{1}\right) .$\emph{\ Set }$\lambda =\lambda \left( \delta \right) =\ln
\delta ^{-1/\left( 4\rho \right) },\beta =\beta \left( \delta \right)
=3\delta ^{\eta /\left( 4\rho \right) }.$ \emph{Let (\ref{5.5}) be true.
Also, assume that the vector function }$W^{\ast }\in B\left( R\right) $\emph{%
. Let }$W_{\min ,\lambda \left( \delta \right) ,\beta \left( \delta \right)
}\in \overline{B\left( R\right) }$\emph{\ be the minimizer of the functional
(\ref{4.9}), which is guaranteed by Theorem 5.3. Also, let the number }$%
\zeta \in \left( 0,\zeta _{0}\right) $\emph{\ in (\ref{5.2}) be the same as
in Theorem 5.3, so as the number }$\omega \in \left( 0,1\right) $\emph{.
Then the following estimates hold:}%
\begin{equation}
\left\Vert W^{\ast }-W_{\min ,\lambda \left( \delta \right) ,\beta \left(
\delta \right) }\right\Vert _{H^{1}\left( \Omega _{\mu +\eta }\right) }\leq
C_{2}\delta ^{\eta /\left( 8\rho \right) },  \label{5.6}
\end{equation}%
\begin{equation}
\left\Vert \Delta W^{\ast }-\Delta W_{n}\right\Vert _{L^{2}\left( \Omega
_{\mu +\eta }\right) }\leq C_{2}\delta ^{\eta /\left( 8\rho \right) },
\label{5.600}
\end{equation}%
\begin{equation}
\left\Vert W^{\ast }-W_{n}\right\Vert _{H^{1}\left( \Omega _{\mu +\eta
}\right) }\leq C_{2}\delta ^{\eta /\left( 8\rho \right) }+\omega
^{n}\left\Vert W_{0}-W_{\min ,\lambda \left( \delta \right) ,\beta \left(
\delta \right) }\right\Vert _{H^{3}\left( \Omega \right) },\text{ }n=1,2,...,
\label{5.60}
\end{equation}%
\begin{equation}
\left\Vert \Delta W^{\ast }-\Delta W_{n}\right\Vert _{L^{2}\left( \Omega
_{\mu +\eta }\right) }\leq C_{2}\delta ^{\eta /\left( 8\rho \right) }+\omega
^{n}\left\Vert W_{0}-W_{\min ,\lambda \left( \delta \right) ,\beta \left(
\delta \right) }\right\Vert _{H^{3}\left( \Omega \right) },\text{ }n=1,2,...
\label{5.61}
\end{equation}%
\emph{In the case of noiseless data with }$\delta =0$\emph{\ one should
replace in (\ref{5.6}), (\ref{5.60}) }$\delta ^{\eta /\left( 8\rho \right) }$%
\emph{\ with }$\sqrt{\beta },$ \emph{where }$\beta =3e^{-\lambda \eta }$%
\emph{\ and }$\lambda \geq \lambda _{2}.$

While (\ref{5.6})-(\ref{5.61}) are convergence estimates for the vector
function $W^{\ast }\left( x\right) ,$ we still need to obtain a convergence
estimate for our target coefficient $a_{0}\left( x\right) $ in equation (\ref%
{eqn:eitv}). This is done in Theorem 5.5. Let $V^{\ast }\left( x\right)
=W^{\ast }\left( x\right) +P^{\ast }\left( x\right) .$ Then $V^{\ast }\left(
x\right) =\left( v_{0}^{\ast }\left( x\right) ,...,v_{N-1}^{\ast }\left(
x\right) \right) ^{T}.$ Let $a_{0}^{\ast }\left( x\right) $ be the exact
coefficient $a_{0}\left( x\right) $ which corresponds to $V^{\ast }\left(
x\right) $ via (\ref{3.90}), i.e.%
\begin{equation}
a_{0}^{\ast }\left( x\right) =-\mathop{\displaystyle \sum }%
\limits_{k=0}^{N-1}\Delta v_{k}^{\ast }\left( x\right) \psi _{k}\left(
s\right) +\left( \mathop{\displaystyle \sum }\limits_{k=0}^{N-1}\nabla
v_{k}^{\ast }\left( x\right) \psi _{k}\left( s\right) \right) ^{2},\text{ }%
x\in \Omega ,s\in \left( 0,1\right) .  \label{5.62}
\end{equation}%
Next, let $V_{\min ,\lambda \left( \delta \right) ,\beta \left( \delta
\right) }\left( x\right) =W_{\min ,\lambda \left( \delta \right) ,\beta
\left( \delta \right) }\left( x\right) +P\left( x\right) =\left( v_{0,\min
,\lambda \left( \delta \right) ,\beta \left( \delta \right) }\left( x\right)
,...,v_{N-1,\min ,\lambda \left( \delta \right) ,\beta \left( \delta \right)
}\left( x\right) \right) ^{T}$ and let%
\begin{equation}
a_{0,\min ,\lambda \left( \delta \right) ,\beta \left( \delta \right)
}\left( x\right) =-\mathop{\displaystyle \sum }\limits_{k=0}^{N-1}\Delta
v_{k,\min ,\lambda \left( \delta \right) ,\beta \left( \delta \right)
}\left( x\right) \psi _{k}\left( \overline{s}\right)  \label{5.63}
\end{equation}%
\[
+\left( \mathop{\displaystyle \sum }\limits_{k=0}^{N-1}\nabla v_{k,\min
,\lambda \left( \delta \right) ,\beta \left( \delta \right) }\left( x\right)
\psi _{k}\left( \overline{s}\right) \right) ^{2},\text{ }x\in \Omega ,%
\overline{s}\in \left( 0,1\right) . 
\]%
Let $V_{n}\left( x\right) =W_{n}\left( x\right) +P\left( x\right) ,$ where
the sequence $\left\{ W_{n}\right\} _{n=0}^{\infty }$ is defined in (\ref%
{5.2}). Then $V_{n}\left( x\right) =\left( v_{0}^{\left( n\right) }\left(
x\right) ,...,v_{N-1}^{\left( n\right) }\left( x\right) \right) ^{T}.$
Define the function $a_{0,n}\left( x\right) $ as 
\begin{equation}
a_{0,n}\left( x\right) =-\mathop{\displaystyle \sum }\limits_{k=0}^{N-1}%
\Delta v_{k}^{\left( n\right) }\left( x\right) \psi _{k}\left( \overline{s}%
\right) +\left( \mathop{\displaystyle \sum }\limits_{k=0}^{N-1}\nabla
v_{k}^{\left( n\right) }\left( x\right) \psi _{k}\left( \overline{s}\right)
\right) ^{2},\text{ }x\in \Omega ,\overline{s}\in \left( 0,1\right) ,
\label{5.64}
\end{equation}%
where $\overline{s}$ is a certain fixed number.

\textbf{Theorem 5.5}. \emph{Assume that conditions of Theorem 5.4 hold. Then
the following analogs of estimates (\ref{5.6})-(\ref{5.61}) are in place:}%
\begin{equation}
\left\Vert a_{0}^{\ast }-a_{0,\min ,\lambda \left( \delta \right) ,\beta
\left( \delta \right) }\right\Vert _{L^{2}\left( \Omega _{\mu +\eta }\right)
}\leq C_{2}\delta ^{\eta /\left( 8\rho \right) },  \label{5.65}
\end{equation}%
\begin{equation}
\left\Vert a_{0}^{\ast }-a_{0,n}\right\Vert _{L^{2}\left( \Omega _{\mu +\eta
}\right) }\leq C_{2}\delta ^{\eta /\left( 8\rho \right) }+\omega
^{n}\left\Vert W_{0}-W_{\min ,\lambda \left( \delta \right) ,\beta \left(
\delta \right) }\right\Vert _{H^{3}\left( \Omega \right) },\text{ }n=1,2,...,
\label{5.66}
\end{equation}%
\emph{where functions }$a_{0}^{\ast },a_{0,\min ,\lambda \left( \delta
\right) ,\beta \left( \delta \right) },a_{0,n}$\emph{\ are defined in (\ref%
{5.62})-(\ref{5.64}).}

\textbf{Remarks 5.1}:

\begin{enumerate}
\item \emph{Theorems 5.4 and 5.5 guarantee that a small neighborhood of the
exact solution is reached if the gradient projection method starts from an
arbitrary point of the ball }$B\left( R\right) $\emph{. Since the radius }$R$%
\emph{\ of this ball is an arbitrary one, then this is global convergence,
see section 1 for our definition of the global convergence. }

\item \emph{The proof of Theorem 5.2 is quite similar with the proof of
theorem 3.1 of \cite{BakKlib}. Theorem 5.3 follows immediately from a
combination of Theorems 5.1 and 5.2 with lemma 2.1 and theorem 2.1 of \cite%
{BakKlib}. Thus, we omit proofs of Theorems 5.2 and 5.3 and focus only on
Theorems 5.1, 5.4 and 5.5. In proofs below }$C_{2}=C_{2}\left( \mu ,\eta
,F,\Psi \left( N\right) ,P,R\right) >0$ \emph{denotes different constants
depending only on listed parameters.}
\end{enumerate}

\subsection{Proof of Theorem 5.1}

\label{sec:5.2}

Let $W_{1},W_{2}\in \overline{B\left( R\right) }$ be two arbitrary points.
Denote $h=W_{2}-W_{1}.$ Hence, $W_{2}=W_{1}+h.$ By the triangle inequality
and (\ref{4.10}) 
\begin{equation}
\left\Vert h\right\Vert _{H^{3}\left( \Omega _{\mu }\right) }\leq 2R.
\label{5.7}
\end{equation}%
We have%
\[
\left[ \Delta W_{1}+\Delta h-F\left( \nabla W_{1}+\nabla P+\nabla h\right)
+\Delta P\right] ^{2}-\left[ \Delta W_{1}-F\left( \nabla W_{1}+\nabla
P\right) +\Delta P\right] ^{2} 
\]%
\begin{equation}
=\left[ \Delta h-\left( F\left( \nabla W_{1}+\nabla P+\nabla h\right)
-F\left( \nabla W_{1}+\nabla P\right) \right) \right]  \label{5.8}
\end{equation}%
\[
\times \left[ \Delta h+2\Delta W_{1}-F\left( \nabla W_{1}+\nabla P+\nabla
h\right) -F\left( \nabla W_{1}+\nabla P\right) +2\Delta P\right] . 
\]%
Recall that the vector function $F\left( \nabla W+\nabla P\right) $ is
quadratic with respect to the derivatives $\partial _{x_{j}}W_{k}\left(
x\right) ,j=1,2,3;k=0,...,N-1.$ Hence, (\ref{5.8}) implies that 
\[
\left[ \Delta W_{2}-F\left( \nabla W_{2}+\nabla P\right) +\Delta P\right]
^{2}-\left[ \Delta W_{1}-F\left( \nabla W_{1}+\nabla P\right) +\Delta P%
\right] ^{2} 
\]%
\begin{equation}
=\Delta h\left[ Q_{1}\left( \nabla W_{1}+\nabla P\right) +2\left( \Delta
W_{1}+\Delta P\right) \right] +\nabla h\left[ Q_{2}\left( \nabla
W_{1}+\nabla P,\Delta W_{1}+\Delta P\right) \right]  \label{5.9}
\end{equation}%
\[
+\left( \Delta h\right) ^{2}+\Delta hD_{1}\left( \nabla W_{1}+\nabla
P,\nabla h\right) +D_{2}\left( \nabla W_{1}+\nabla P,\nabla h\right) . 
\]%
In (\ref{5.9}) vector functions $Q_{1},Q_{2},D_{1},D_{2}$ are continuous
with respect to their indicated variables. In addition, (\ref{4.8}) and (\ref%
{5.7}) imply that the following estimates are valid for vector functions $%
D_{1}\left( \nabla W_{1}+\nabla P,\nabla h\right) ,D_{2}\left( \nabla
W_{1}+\nabla P,\nabla h\right) $:%
\begin{equation}
\left\vert D_{1}\left( \nabla W_{1}+\nabla P,\nabla h\right) \right\vert
\leq C_{2}\left( \left\vert \nabla h\right\vert +\left\vert \nabla
h\right\vert ^{2}\right) ,  \label{5.10}
\end{equation}%
\begin{equation}
\left\vert D_{2}\left( \nabla W_{1}+\nabla P,\nabla h\right) \right\vert
\leq C_{2}\left\vert \nabla h\right\vert ^{2},\text{ }j=1,2.  \label{5.100}
\end{equation}%
In the second line of (\ref{5.9}), we single out the part which is linear
with respect to $h$. On the other hand, using (\ref{4.8}), (\ref{5.10}), (%
\ref{5.100}) and Cauchy-Schwarz inequality, we obtain the following estimate
from the below for the expression in the third line of (\ref{5.9}):%
\begin{equation}
\left( \Delta h\right) ^{2}+\Delta hD_{1}\left( \nabla W_{1}+\nabla P,\nabla
h\right) +D_{2}\left( \nabla W_{1}+\nabla P,\nabla h\right) \geq \frac{1}{2}%
\left( \Delta h\right) ^{2}-C_{2}\left( \nabla h\right) ^{2}.  \label{5.11}
\end{equation}%
In addition, the following estimate from the above follows from (\ref{4.8}),
(\ref{5.9}), (\ref{5.10}) and (\ref{5.100}):%
\begin{equation}
\left\vert \left( \Delta h\right) ^{2}+\Delta hD_{1}\left( \nabla
W_{1}+\nabla P,\nabla h\right) +D_{2}\left( \nabla W_{1}+\nabla P,\nabla
h\right) \right\vert  \label{5.12}
\end{equation}%
\[
\leq C_{2}\left[ \left( \Delta h\right) ^{2}+\left( \nabla h\right) ^{2}%
\right] . 
\]

Thus, (\ref{4.9}) and (\ref{5.9}) imply that 
\[
J_{\lambda ,\beta }\left( W_{1}+h\right) -J_{\lambda ,\beta }\left(
W_{1}\right) 
\]%
\[
=e^{-2\lambda \left( \mu +\eta \right) }\mathop{\displaystyle \int}%
\limits_{\Omega _{\mu }}\left\{ \Delta h\left[ Q_{1}\left( \nabla
W_{1}+\nabla P\right) +2\Delta P\right] +\nabla h\left[ Q_{2}\left( \nabla
W_{1}+\nabla P,\Delta P\right) \right] \right\} e^{2\lambda r}dx 
\]%
\begin{equation}
+2\beta \left[ h,W_{1}\right]  \label{5.13}
\end{equation}%
\[
+e^{-2\lambda \left( \mu +\eta \right) }\mathop{\displaystyle \int}%
\limits_{\Omega _{\mu }}\left[ \left( \Delta h\right) ^{2}+\Delta
hD_{1}\left( \nabla W_{1}+\nabla P,\nabla h\right) +D_{2}\left( \nabla
W_{1}+\nabla P,\nabla h\right) \right] e^{2\lambda r}dx 
\]%
\[
+\beta \left\Vert h\right\Vert _{H^{3}\left( \Omega _{\mu }\right) }^{2}. 
\]%
The expression in the second line of (\ref{5.13}) is generated by the second
line of (\ref{5.9}), and it is linear with respect to $h$. Actually, the sum
of the second and third lines of (\ref{5.13}) is a linear functional with
respect to $h,$ and we denote it $Lin\left( W_{1}\right) \left( h\right) .$%
In addition, the following estimate holds 
\[
\left\vert Lin\left( W_{1}\right) \left( h\right) \right\vert \leq C_{2}\exp
\left( 2\lambda \left( \rho -\mu -\eta \right) \right) \left\Vert
h\right\Vert _{H^{3}\left( \Omega _{\mu }\right) }. 
\]%
Hence, $Lin\left( W_{1}\right) \left( h\right) :H_{0}^{3}\left( \Omega _{\mu
}\right) \rightarrow \mathbb{R}$ is a bounded linear functional with respect
to $h$. Hence, by Riesz theorem there exists a vector function $Y\left(
x\right) \in H_{0}^{3}\left( \Omega _{\mu }\right) $ such that 
\begin{equation}
Lin\left( W_{1}\right) \left( h\right) =\left[ Y,h\right] .  \label{5.14}
\end{equation}%
Also, it follows from (\ref{5.10}) and (\ref{5.13}) that if $\left\Vert
h\right\Vert _{H^{3}\left( \Omega _{\mu }\right) }<1,$ then the following
estimate holds 
\begin{equation}
\left\vert J_{\lambda ,\beta }\left( W_{1}+h\right) -J_{\lambda ,\beta
}\left( W_{1}\right) -Lin\left( W_{1}\right) \left( h\right) \right\vert
\leq C_{2}\exp \left( 2\lambda \left( \rho -\mu -\eta \right) \right)
\left\Vert h\right\Vert _{H^{3}\left( \Omega _{\mu }\right) }^{2}.
\label{5.15}
\end{equation}%
Thus, using (\ref{5.14}) and (\ref{5.15}), we obtain that the Frech\'{e}t
derivative $J_{\lambda ,\beta }^{\prime }\left( W_{1}\right) $ of the
functional $J_{\lambda ,\beta }\left( W\right) $ exists at the point $W_{1}$
and $J_{\lambda ,\beta }^{\prime }\left( W_{1}\right) =Y\left( x\right) .$
Even though the existence of the Frech\'{e}t derivative $J_{\lambda ,\beta
}^{\prime }\left( W_{1}\right) $ is proved here only for the case when $%
W_{1} $ is an interior point of the ball $B\left( R\right) ,$ still since $%
R>0$ is an arbitrary number, then actually this existence is proved for an
arbitrary point $W_{1}\in H_{0}^{3}\left( \Omega _{\mu }\right) .$

We now need to prove the strict convexity estimate (\ref{5.1}). To do this,
we will use the Carleman estimate of Theorem 3.1. Using (\ref{5.11}) and (%
\ref{5.13}), we obtain 
\begin{equation}
J_{\lambda ,\beta }\left( W_{1}+h\right) -J_{\lambda ,\beta }\left(
W_{1}\right) -J_{\lambda ,\beta }^{\prime }\left( W_{1}\right) \left(
h\right)  \label{5.16}
\end{equation}%
\[
\geq \frac{1}{2}e^{-2\lambda \left( \mu +\eta \right) }%
\mathop{\displaystyle
\int}\limits_{\Omega _{\mu }}\left( \Delta h\right) ^{2}e^{2\lambda
r}dx-C_{2}e^{-2\lambda \left( \mu +\eta \right) }\mathop{\displaystyle \int}%
\limits_{\Omega _{\mu }}\left( \nabla h\right) ^{2}e^{2\lambda r}dx+\beta
\left\Vert h\right\Vert _{H^{3}\left( \Omega _{\mu }\right) }^{2}. 
\]%
Next, using (\ref{3.122}), we obtain from (\ref{5.16}) 
\[
J_{\lambda ,\beta }\left( W_{1}+h\right) -J_{\lambda ,\beta }\left(
W_{1}\right) -J_{\lambda ,\beta }^{\prime }\left( W_{1}\right) \left(
h\right) 
\]%
\begin{equation}
\geq \frac{1}{4}e^{-2\lambda \left( \mu +\eta \right) }\mathop{\displaystyle
\int}\limits_{\Omega _{\mu }}\left( \Delta h\right) ^{2}e^{2\lambda
r}dx+C\lambda e^{-2\lambda \left( \mu +\eta \right) }\mathop{\displaystyle
\int}\limits_{\Omega _{\mu }}\left( \nabla h\right) ^{2}e^{2\lambda r}dx
\label{5.18}
\end{equation}%
\[
-C_{2}e^{-2\lambda \left( \mu +\eta \right) }\mathop{\displaystyle \int}%
\limits_{\Omega _{\mu }}\left( \nabla h\right) ^{2}e^{2\lambda r}dx+C\lambda
^{3}e^{-2\lambda \left( \mu +\eta \right) }\mathop{\displaystyle
\int}\limits_{\Omega _{\mu }}h^{2}e^{2\lambda r}dx 
\]%
\[
-C\lambda ^{3}e^{-2\lambda \eta }\left\Vert h\right\Vert _{H^{2}\left(
\Omega _{\mu }\right) }^{2}+\beta \left\Vert h\right\Vert _{H^{3}\left(
\Omega _{\mu }\right) }^{2}. 
\]%
Choose $\lambda _{2}=\lambda _{2}\left( \mu ,\eta ,F,N,P,R\right) \geq
\lambda _{0}>0$ so large that $C\lambda _{2}>2C_{2}$ and also that $C\lambda
^{3}e^{-2\lambda \eta }<e^{-\lambda \eta },\forall \lambda \geq \lambda
_{2}. $ Recalling (\ref{5.0}) and using $\Omega _{\mu +\eta }\subset \Omega
_{\mu },$ we obtain from (\ref{5.18})%
\[
J_{\lambda ,\beta }\left( W_{1}+h\right) -J_{\lambda ,\beta }\left(
W_{1}\right) -J_{\lambda ,\beta }^{\prime }\left( W_{1}\right) \left(
h\right) \geq \frac{1}{4}e^{-2\lambda \left( \mu +\eta \right) }%
\mathop{\displaystyle \int}\limits_{\Omega _{\mu }}\left( \Delta h\right)
^{2}e^{2\lambda r}dx 
\]%
\begin{equation}
+C_{2}e^{-2\lambda \left( \mu +\eta \right) }\mathop{\displaystyle \int}%
\limits_{\Omega _{\mu }}\left( \nabla h\right) ^{2}e^{2\lambda r}dx+C\lambda
^{3}e^{-2\lambda \left( \mu +\eta \right) }\mathop{\displaystyle \int}%
\limits_{\Omega _{\mu }}h^{2}e^{2\lambda r}dx+\frac{\beta }{2}\left\Vert
h\right\Vert _{H^{3}\left( \Omega _{\mu }\right) }^{2}  \label{5.19}
\end{equation}%
\[
\geq \frac{1}{4}e^{-2\lambda \left( \mu +\eta \right) }%
\mathop{\displaystyle
\int}\limits_{\Omega _{\mu }}\left( \Delta h\right) ^{2}e^{2\lambda
r}dx+C_{2}e^{-2\lambda \left( \mu +\eta \right) }\mathop{\displaystyle \int}%
\limits_{\Omega _{\mu +\eta }}\left[ \left( \nabla h\right) ^{2}+h^{2}\right]
e^{2\lambda r}dx+\frac{\beta }{2}\left\Vert h\right\Vert _{H^{3}\left(
\Omega _{\mu }\right) }^{2}. 
\]%
Next, $e^{2\lambda r}\geq e^{2\lambda \left( \mu +\eta \right) }$ for $x\in
\Omega _{\mu +\eta }.$ Hence, 
\begin{equation}
\frac{1}{4}e^{-2\lambda \left( \mu +\eta \right) }\mathop{\displaystyle \int}%
\limits_{\Omega _{\mu }}\left( \Delta h\right) ^{2}e^{2\lambda
r}dx+e^{-2\lambda \left( \mu +\eta \right) }\mathop{\displaystyle \int}%
\limits_{\Omega _{\mu +\eta }}\left[ \left( \nabla h\right) ^{2}+h^{2}\right]
e^{2\lambda r}dx\geq  \label{5.20}
\end{equation}%
\[
\frac{1}{4}\mathop{\displaystyle \int}\limits_{\Omega _{\mu +\eta }}\left(
\Delta h\right) ^{2}dx+\mathop{\displaystyle \int}\limits_{\Omega _{\mu
+\eta }}\left[ \left( \nabla h\right) ^{2}+h^{2}\right] dx. 
\]%
Thus, (\ref{5.19}) and (\ref{5.20}) imply that%
\begin{equation}
J_{\lambda ,\beta }\left( W_{1}+h\right) -J_{\lambda ,\beta }\left(
W_{1}\right) -J_{\lambda ,\beta }^{\prime }\left( W_{1}\right) \left(
h\right) \geq  \label{5.300}
\end{equation}%
\[
\frac{1}{4}\mathop{\displaystyle \int}\limits_{\Omega _{\mu +\eta }}\left(
\Delta h\right) ^{2}dx+C_{2}\mathop{\displaystyle \int}\limits_{\Omega _{\mu
+\eta }}\left[ \left( \nabla h\right) ^{2}+h^{2}\right] dx+\frac{\beta }{2}%
\left\Vert h\right\Vert _{H^{3}\left( \Omega _{\mu }\right) }^{2}. 
\]%
$\square $

\subsection{Proof of Theorem 5.4}

\label{sec:5.3}

Temporary change notation for the functional (\ref{4.9}) as 
\begin{equation}
J_{\lambda \left( \delta \right) ,\beta \left( \delta \right) }\left(
W+P\right) =e^{-2\lambda \left( \mu +\eta \right) }\mathop{\displaystyle
\int}\limits_{\Omega _{\mu }}\left[ \Delta W+\Delta P-F\left( \nabla
W+\nabla P\right) \right] ^{2}e^{2\lambda r}dx  \label{5.21}
\end{equation}%
\[
+\beta \left\Vert W+P\right\Vert _{H^{3}\left( \Omega _{\mu }\right) }^{2}. 
\]%
Obviously 
\[
e^{-2\lambda \left( \mu +\eta \right) }\mathop{\displaystyle \int}%
\limits_{\Omega _{\mu }}\left[ \Delta W^{\ast }+\Delta P^{\ast }-F\left(
\nabla W^{\ast }+\nabla P^{\ast }\right) \right] ^{2}e^{2\lambda r}dx=0. 
\]%
Hence, by (\ref{5.5}) and (\ref{5.21}) 
\begin{equation}
J_{\lambda \left( \delta \right) ,\beta \left( \delta \right) }\left(
W^{\ast }+P^{\ast }\right) =\beta \left( \delta \right) \left\Vert W^{\ast
}+P^{\ast }\right\Vert _{H^{3}\left( \Omega _{\mu }\right) }^{2}\leq
C_{2}\beta \left( \delta \right) .  \label{5.23}
\end{equation}%
By (\ref{5.21})%
\[
J_{\lambda \left( \delta \right) ,\beta \left( \delta \right) }\left(
W^{\ast }+P\right) -J_{\lambda \left( \delta \right) ,\beta \left( \delta
\right) }\left( W^{\ast }+P^{\ast }\right) = 
\]%
\begin{equation}
e^{-2\lambda \left( \mu +\eta \right) }\mathop{\displaystyle \int}%
\limits_{\Omega _{\mu }}\left[ \Delta W^{\ast }+\Delta P^{\ast }+\Delta 
\widetilde{P}-F\left( \nabla W^{\ast }+\nabla P^{\ast }+\nabla \widetilde{P}%
\right) \right] ^{2}e^{2\lambda r}dx  \label{5.25}
\end{equation}%
\[
-e^{-2\lambda \left( \mu +\eta \right) }\mathop{\displaystyle \int}%
\limits_{\Omega _{\mu }}\left[ \Delta W^{\ast }+\Delta P^{\ast }-F\left(
\nabla W^{\ast }+\nabla P^{\ast }\right) \right] ^{2}e^{2\lambda r}dx 
\]%
\[
+\beta \left[ \widetilde{P},2W^{\ast }+P+P^{\ast }\right] . 
\]%
Recall that $F\left( \nabla V\right) $ is a quadratic vector function with
respect to the derivatives $\partial _{x_{j}}v_{k}\left( x\right) .$ Hence, (%
\ref{5.5}) and (\ref{5.25}) imply that%
\begin{equation}
\left\vert J_{\lambda \left( \delta \right) ,\beta \left( \delta \right)
}\left( W^{\ast }+P\right) -J_{\lambda \left( \delta \right) ,\beta \left(
\delta \right) }\left( W^{\ast }+P^{\ast }\right) \right\vert \leq
C_{2}\delta e^{2\lambda \rho }+C_{2}\delta \beta .  \label{5.26}
\end{equation}%
Next, 
\[
\left\vert J_{\lambda \left( \delta \right) ,\beta \left( \delta \right)
}\left( W^{\ast }+P\right) -J_{\lambda \left( \delta \right) ,\beta \left(
\delta \right) }\left( W^{\ast }+P^{\ast }\right) \right\vert \geq
J_{\lambda \left( \delta \right) ,\beta \left( \delta \right) }\left(
W^{\ast }+P\right) -J_{\lambda \left( \delta \right) ,\beta \left( \delta
\right) }\left( W^{\ast },P^{\ast }\right) . 
\]%
Hence, using (\ref{5.23}) and (\ref{5.26}) and keeping in mind that $%
C_{2}\delta \beta <C_{2}\beta ,$ we obtain 
\begin{equation}
J_{\lambda \left( \delta \right) ,\beta \left( \delta \right) }\left(
W^{\ast }+P\right) \leq C_{2}\delta e^{2\lambda \rho }+C_{2}\beta .
\label{5.27}
\end{equation}%
Since $\lambda \left( \delta \right) =\ln \delta ^{-1/\left( 4\rho \right) }$
and $\delta <\delta _{1}<\min \left( e^{-4\rho \lambda _{2}},3^{-4\rho /\eta
}\right) ,$ then $\lambda \left( \delta \right) >\lambda _{2}$ and also $%
\delta e^{2\lambda \rho }=\sqrt{\delta }.$ Next, since $\beta =3\delta
^{\eta /\left( 4\rho \right) },$ then condition (\ref{5.0}) is fulfilled.
Also, since $3\delta ^{\eta /\left( 4\rho \right) }>\sqrt{\delta },$ then \ $%
\delta e^{2\lambda \rho }+\beta \leq 2\delta ^{\eta /\left( 4\rho \right) }.$

Hence, using (\ref{5.27}), we obtain 
\begin{equation}
J_{\lambda ,\beta }\left( W^{\ast }+P\right) \leq C_{2}\delta ^{\eta /\left(
4\rho \right) }.  \label{5.28}
\end{equation}%
Since by (\ref{5.200}) $\left[ J_{\lambda \left( \delta \right) ,\beta
\left( \delta \right) }^{\prime }\left( W_{\min ,\lambda \left( \delta
\right) ,\beta \left( \delta \right) }\right) ,W^{\ast }-W_{\min ,\lambda
\left( \delta \right) ,\beta \left( \delta \right) }\right] \geq 0,$ then,
using (\ref{5.28}), we obtain 
\[
J_{\lambda \left( \delta \right) ,\beta \left( \delta \right) }\left(
W^{\ast }\right) -J_{\lambda \left( \delta \right) ,\beta \left( \delta
\right) }\left( W_{\min ,\lambda \left( \delta \right) ,\beta \left( \delta
\right) }\right) -J_{\lambda \left( \delta \right) ,\beta \left( \delta
\right) }^{\prime }\left( W_{\min ,\lambda \left( \delta \right) ,\beta
\left( \delta \right) }\right) \left( W^{\ast }-W_{\min ,\lambda }\right)
\leq C_{2}\delta ^{\eta /\left( 4\rho \right) }. 
\]%
Hence, by (\ref{5.1}) 
\[
\left\Vert \Delta W^{\ast }-\Delta W_{\min ,\lambda \left( \delta \right)
,\beta \left( \delta \right) }\right\Vert _{L^{2}\left( \Omega _{\mu +\eta
}\right) }^{2}+\left\Vert W^{\ast }-W_{\min ,\lambda \left( \delta \right)
,\beta \left( \delta \right) }\right\Vert _{H^{1}\left( \Omega _{\mu +\eta
}\right) }^{2}\leq C_{2}\delta ^{\eta /\left( 4\rho \right) }, 
\]%
from which (\ref{5.6}) and (\ref{5.600}) follow.

We now prove (\ref{5.60}). Using triangle inequality (\ref{5.201}) and (\ref%
{5.6}), we obtain for $n\geq 1$ 
\[
\left\Vert W^{\ast }-W_{n}\right\Vert _{H^{1}\left( \Omega _{\mu +\eta
}\right) }\leq \left\Vert W^{\ast }-W_{\min ,\lambda \left( \delta \right)
,\beta \left( \delta \right) }\right\Vert _{H^{1}\left( \Omega _{\mu +\eta
}\right) }+\left\Vert W_{\min ,\lambda \left( \delta \right) ,\beta \left(
\delta \right) }-W_{n}\right\Vert _{H^{1}\left( \Omega _{\mu +\eta }\right)
} 
\]%
\[
\leq C_{2}\delta ^{\eta /\left( 8\rho \right) }+\omega ^{n}\left\Vert
W_{0}-W_{\min ,\lambda \left( \delta \right) ,\beta \left( \delta \right)
}\right\Vert _{H^{3}\left( \Omega \right) }, 
\]%
which proves (\ref{5.60}). The proof of (\ref{5.61}) is completely similar. $%
\square $

\subsection{Proof of Theorem 5.5}

\label{sec:5.4}

Subtracting (\ref{5.63}) from (\ref{5.62}), we obtain%
\begin{equation}
\left\vert a_{0}^{\ast }\left( x\right) -a_{0,\min ,\lambda \left( \delta
\right) ,\beta \left( \delta \right) }\left( x\right) \right\vert \leq C_{2}%
\mathop{\displaystyle \sum }\limits_{k=0}^{N-1}\left\vert \Delta v_{k}^{\ast
}\left( x\right) -\Delta v_{k,\min ,\lambda \left( \delta \right) ,\beta
\left( \delta \right) }\left( x\right) \right\vert  \label{5.29}
\end{equation}%
\[
+C_{2}\mathop{\displaystyle \sum }\limits_{k=0}^{N-1}\left\vert \nabla
v_{k}^{\ast }\left( x\right) -\nabla v_{k,\min ,\lambda \left( \delta
\right) ,\beta \left( \delta \right) }\left( x\right) \right\vert \left\vert
\nabla v_{k}^{\ast }\left( x\right) +\nabla v_{k,\min ,\lambda \left( \delta
\right) ,\beta \left( \delta \right) }\left( x\right) \right\vert . 
\]%
Since vector functions $W_{\min ,\lambda \left( \delta \right) ,\beta \left(
\delta \right) },W^{\ast }\in \overline{B\left( R\right) }$, then (\ref{4.10}%
) implies that $\left\vert \nabla v_{k}^{\ast }\left( x\right) +\nabla
v_{k,\min ,\lambda \left( \delta \right) ,\beta \left( \delta \right)
}\left( x\right) \right\vert \leq C_{2}.$ Hence, (\ref{5.65}) follows from (%
\ref{5.6}), (\ref{5.600}) and (\ref{5.29}). The proof of (\ref{5.66}) is
completely similar. $\square $

\section{Numerical studies}

\label{sec:6}

We have applied the above technique to numerical studies of the inverse EIT
problem in the 2D case. Recall that even though theorems 5.1-5.4 are
formulated only in the 3D case, their direct analogs are also valid in the
2D case due to the Carleman estimate of Theorem 3.2, see beginning of
section 4. In this section we describe our numerical results. Hence, in this
section%
\[
\Omega =\left\{ r\in \left( 0,\rho \right) \right\} \subset \mathbb{R}%
^{2},\Omega _{\mu }=\left\{ r\in \left( \mu ,\rho \right) \right\} \subset
\Omega . 
\]%
We have found in our computations that the influence of the regularization
parameter $\beta $ in (\ref{4.9}) is not essential. Hence, we set $\beta :=0$
in our computational examples.

\subsection{Some details of the numerical implementation}

\label{sec:6.1}

In all our numerical examples 
\[
G=\{x_{1}^{2}+x_{2}^{2}<5\},\Omega =\{x_{1}^{2}+x_{2}^{2}\leq 1\}\text{ and }%
\Omega _{\mu }=\left\{ r\in \left( 0.01,1\right) \right\} \subset \Omega . 
\]%
We measure the data on the whole boundary $\partial \Omega =S_{1}$. The
source runs over the circle $C^{\left( s\right) }=\{x_{1}^{2}+x_{2}^{2}=4\}$%
. In other words, in polar coordinates 
\begin{equation}
x_{s}=\left( r,s\right) =\left( 2,s\right) ,s=\varphi \in \left( 0,2\pi
\right) ,x_{s}\in C^{\left( s\right) }.  \label{6.20}
\end{equation}

\begin{figure}[tbp]
\begin{center}
\includegraphics[width=6cm]{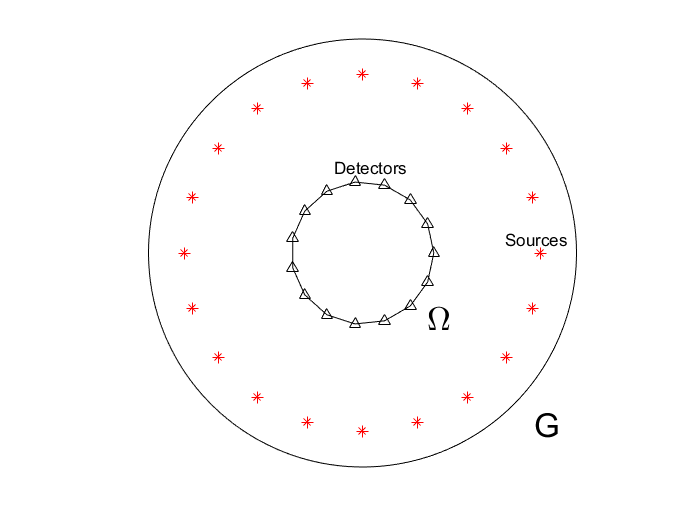}
\end{center}
\caption{A schematic diagram of
domains $G,\Omega ,$ sources and detectors.}
\label{domain}
\end{figure}

However, when constructing the required orthonormal basis $\left\{ \psi
_{n}\left( s\right) \right\} _{n=0}^{\infty },$ we still have used functions 
$\left\{ s^{n}e^{s}\right\} _{n=0}^{\infty },$ i.e. we did not impose the
periodicity condition on this basis. The source function $f\left( x\right) $
in our case is the bump function below:%
\[
f\left( x-x^{\left( s\right) }\right) =\left\{ 
\begin{array}{l}
\frac{1}{\varepsilon }\exp \left( -\frac{1}{1-\left\vert x-x_{s}\right\vert
^{2}/\varepsilon }\right) ,\text{ }if\text{ }\left( x-x_{s}\right)
^{2}<\varepsilon , \\ 
0,\text{ otherwise.}%
\end{array}%
\right. 
\]%
We have chosen $\varepsilon =0.01$.

We use 32 sources and 32 detectors. In examples 1-5 both sources and
detectors are uniformly distributed over the whole circle $%
\{x_{1}^{2}+x_{2}^{2}=4\}$ and the whole circle $S_{1}=%
\{x_{1}^{2}+x_{2}^{2}=1\}$ respectively. However, this changes in Example 6
(see below).

To solve the forward problem (\ref{eqn:eitres}), we have used the standard
FEM. However, to minimize functional (\ref{4.9}), we have written the
differential operators in it via finite differences. Thus, we have not
committed \textquotedblleft inverse crime". To use the finite differences,
we have discretized the domain $\Omega _{\mu }$ in polar coordinates using
the uniform finite difference mesh$.$ Next, we have used the gradient
descent method to minimize functional (\ref{4.9}) with respect to the values
of the vector function $W\left( r,\varphi \right) $ at grid points. As the
basis $\psi _{k}$ is not periodic over $[0,2\pi ]$, we treat numerically $%
s=0 $ and $s=2\pi $ as two different discrete points.

As to the choice of the parameter $\lambda ,$even though the above theory
works only for sufficiently large values of $\lambda $, we have established
in our computational experiments that the choice 
\begin{equation}
\lambda =1  \label{6.3}
\end{equation}%
is sufficient for all six tests we have performed. We have also tested three
different values of the number $N$ terms in the series (\ref{3.3}): 
\begin{equation}
N=4,6,8.  \label{6.4}
\end{equation}%
Our computational results indicate that $N=8$ is the best choice out of
these three.

\textbf{Remark 6.1.} \emph{The choice (\ref{6.3}) of the parameter }$\lambda 
$\emph{\ corresponds well with the observations of previous publications on
numerical studies of the convexification method, both for coefficient
inverse problems with the single location of the source \cite%
{KlibThanh,KlibKol1,KlibKol2} and for ill-posed problems for quasilinear
parabolic equations \cite{BakKlib,KlibYag}. This observation is that not
large values of }$\lambda $\emph{\ can be chosen in computations. }

\subsection{A multi-level method of the minimization of functional (\protect
\ref{4.9})}

\label{sec:6.2}

We have found in our computational experiments that the gradient descent
method for our weighted Tikhonov-like functional (\ref{4.9}) converges
rapidly on a coarse mesh. This provides us with a rough image. Hence, we
have implemented a multi-level method \cite{Zouli2007}. Let $%
M_{h_{1}}\subset M_{h_{2}}...\subset M_{h_{K}}$ be nested finite difference
meshes, i.e. $M_{h_{k}}$ is a refinement of $M_{h_{k-1}}$ for $k\leq K$. Let 
$P_{h_{k}}$ be the corresponding finite difference functional space. One the
first level $M_{h_{1}}$, we solve the discrete optimization problem. In
other words, let $V_{h_{1},\min }$ be the minimizer of the following
functional which is found via the gradient descent method%
\begin{equation}
J_{\lambda }^{\left( h_{1}\right) }\left( W_{h_{1}}\right) =e^{-2\lambda
\left( \mu +\eta \right) }\mathop{\displaystyle \int}\limits_{\Omega _{\mu }}%
\left[ \Delta W_{h_{1}}+\Delta P-F\left( \nabla W_{h_{1}}+\nabla P\right) %
\right] ^{2}e^{2\lambda r}dx,  \label{6.5}
\end{equation}%
where the integral is understood in the discrete sense. Then we interpolate
the minimizer $W_{h_{1},\min ,\lambda }$ on the finer mesh $M_{h_{2}}$ and
take the resulting vector function $W_{h_{2},\text{int}}$ as the starting
point of the gradient descent method of the optimization of the direct
analog of functional (\ref{6.5}) in which $h_{1}$ is replaced with $h_{2}$
and $W_{h_{1}}$ is replaced with $W_{h_{2}}.$This process was repeated until
we got the minimizer $W_{h_{K},\min }$ on the $K_{th}$ level on the mesh $%
M_{h_{K}}$.

Since $(r,\varphi )\in (0,1)\times (0,2\pi )$, then our first level $%
M_{h_{1}}$ is set to be the uniform mesh with the mesh size in the $r$
direction to be $1/4$ and the mesh size in the $\varphi $ direction to be $%
2\pi /8$. For each mesh refinement, we will refine the mesh in both $r$
direction and $\varphi $ direction in a way that we set the mesh size of the refined mesh in both direction to be 1/2 of the previous mesh sizes. On each level $M_{h_{k}}$, as soon as we see that $\Vert \nabla
J_{\lambda }^{(h_{k})}(W_{h_{k}})\Vert <2\times 10^{-2}$, we refine the mesh
and compute the solution on the next level $M_{h_{k+1}}$. In the end, we
compute $a_{0}(x)$ using the relation \eqref{3.90} with $s=0$.

Our starting point $W^{\left( 0\right) }\left( r,\varphi \right) $ for the
vector function $W\left( r,\varphi \right) $ for the gradient descent method
on the coarse mesh $M_{h_{1}}$ is set to be the background solution $%
W^{\left( 0\right) }(r,\varphi ,1)$ which corresponds to the solution of the
problem (\ref{eqn:eitres}) with $\sigma \left( x\right) \equiv 1$. Hence,
our starting point is not located in a small neighborhood of the exact
solution.

\subsection{Numerical testing}

\label{sec:6.3}

In the tests of this section, we demonstrate the efficiency of our numerical
method for imaging of small inclusions as well as for imaging of a smoothly
varying function $\sigma \left( x\right) ,$ i.e. a \textquotedblleft
stretched" inclusion with a wide range of change of the conductivity inside
of it. In particular, we test the case of a rather high contrast 5:1 of the
inclusion. In all tests the background value of the conductivity is $\sigma
_{bkgr}=1.$ In addition, we test the influence of the number $N$ in (\ref%
{6.4}). We also test the effects of both: the data given only on a part of
the boundary and the source running only along a part of the circle $\left\{
r=2\right\} .$ In Tests 1-6 we have stopped on the 3rd mesh refinement for
all three values of $N$ listed in (\ref{6.4}) (except for test 4 where $N=8$). The reason of stopping on the 3rd mesh refinement is that images were changing very insignificantly when on the 3rd mesh refinement, as compared with the second. 

All necessary derivatives of the data were calculated using finite
differences, just as in previous above cited publications of the first
author with coauthors about the convexification \cite{KlibThanh,KlibKol2}
with numerical results in them, including the one with noisy experimental
data \cite{KlibKol1}. Just as in those works, we have not observed
instabilities due to the differentiation, most likely because the step sizes
of finite differences were not too small.

\begin{figure}[tbp]
\begin{center}
\begin{tabular}{cc}
\includegraphics[width=5cm]{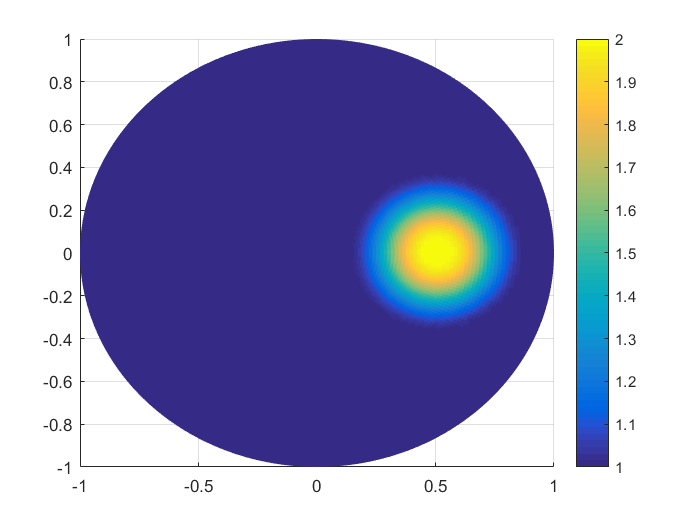} & %
\includegraphics[width=5cm]{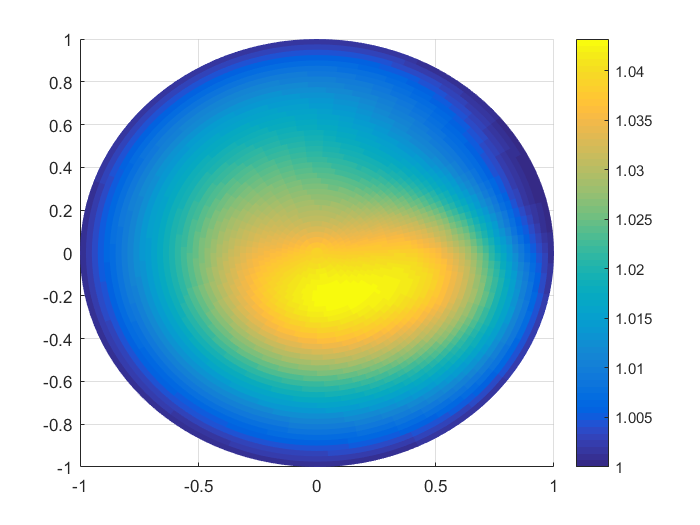} \\ 
(a) True $\sigma$ & (b) $N=4$ \\ 
\includegraphics[width=5cm]{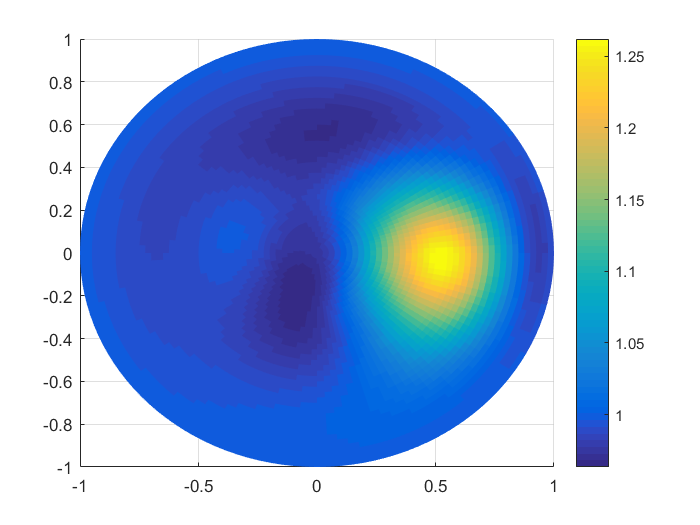} & %
\includegraphics[width=5cm]{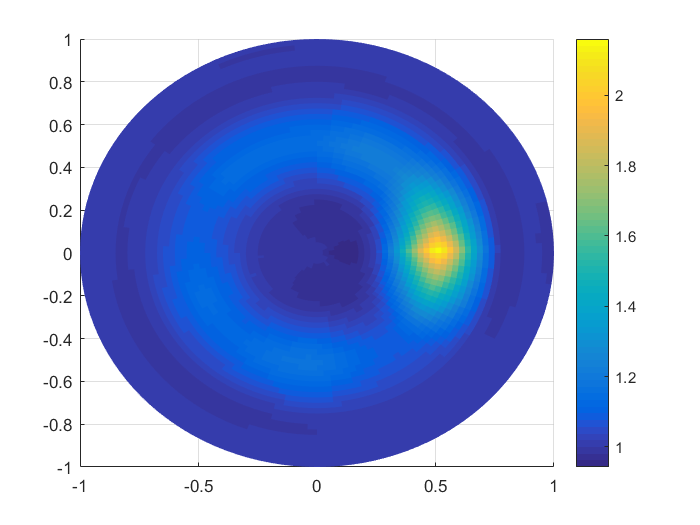} \\ 
(c) $N=6$ & (d) $N=8$%
\end{tabular}%
\end{center}
\caption{\emph{Results of Test 1.\ Imaging of one inclusion with } $\protect%
\sigma =2$\emph{\ in it and }$\protect\sigma =1$ \emph{outside. Hence, the
inclusion/background contrast is 2:1. We have stopped at the 3rd mesh
refinement for all three values of }$N$\emph{\ listed in (\protect\ref{6.4}%
). a) Correct image. b) Computed image for }$N=4$\emph{. c) Computed image
for }$N=6$\emph{. d) Computed image for }$N=8$\emph{. Both the correct
contrast and correct location are achieved at }$N=8$\emph{. } }
\label{example1}
\end{figure}

\textbf{Test 1}. First, we test the reconstruction by our method of a single
inclusion depicted on Figure \ref{example1} a). $\sigma =2$ inside of this
inclusion and $\sigma =1$ outside. Hence, the inclusion/background contrast
is 2:1. The best result is achieved at $N=8$, see Figures \ref{example1}.

\begin{figure}[tbp]
\begin{center}
\begin{tabular}{cc}
\includegraphics[width=5cm]{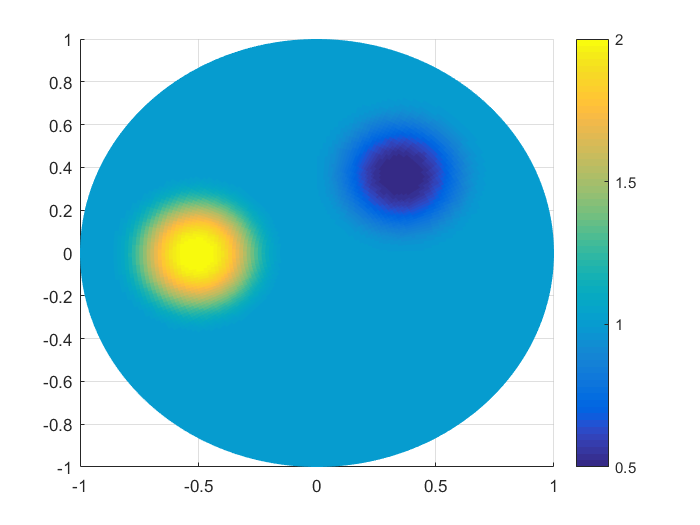} & %
\includegraphics[width=5cm]{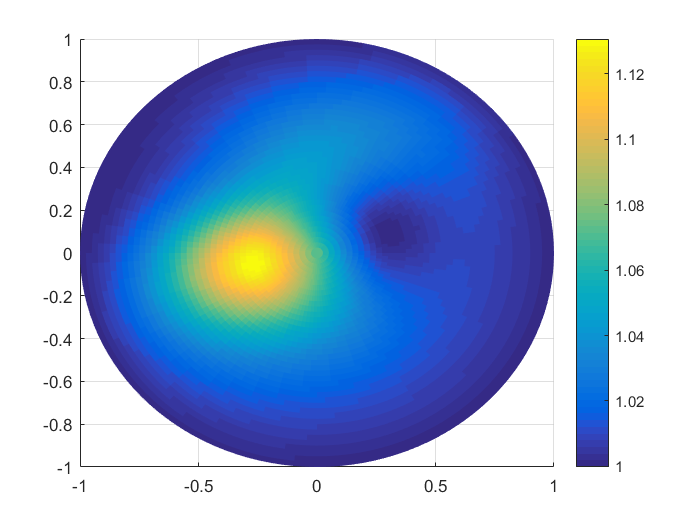} \\ 
(a) True $\sigma$ & (b) $N=4$ \\ 
\includegraphics[width=5cm]{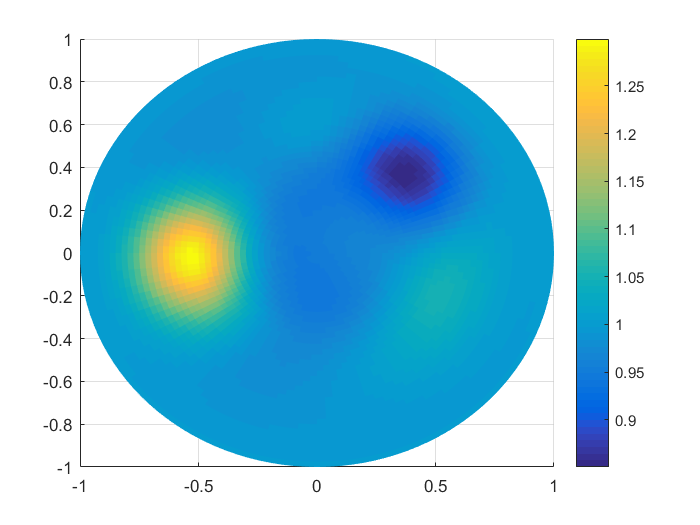} & %
\includegraphics[width=5cm]{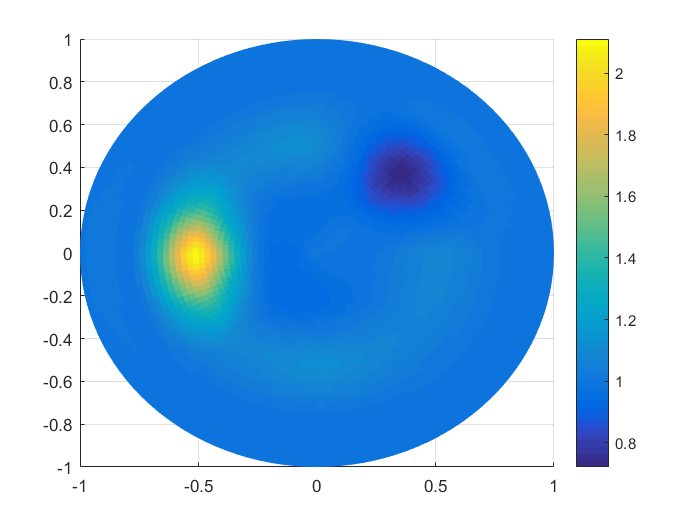} \\ 
(c) $N=6$ & (d) $N=8$%
\end{tabular}%
\end{center}
\caption{\emph{Results of Test 2.\ Imaging of two inclusions. Here, }$%
\protect\sigma =2$\emph{\ in the left inclusion, }$\protect\sigma =0.5$ 
\emph{in the right inclusion and }$\protect\sigma =1$ \emph{otherwise.
Hence, the inclusion/background contrast is 2:1 in the left inclusion and is
0.5:1 in the right inclusion. This means that the electric conductivity of
the left inclusion is higher than the one of the background and it is lower
of the right inclusion. We have stopped on the 3rd mesh refinement for all
three values of }$N$\emph{\ listed in (\protect\ref{6.4}). a) Correct image.
b) Computed image for }$N=4$\emph{. c) Computed image for }$N=6$\emph{. d)
Computed image for }$N=8$\emph{, which is the best one out of three. }}
\label{example2}
\end{figure}

\textbf{Test 2}. We test now the performance of our method for imaging of
two inclusions depicted on Figure \ref{example2} a). $\sigma =2$ inside of
each inclusion and $\sigma =1$ outside of these inclusions. See Figures \ref%
{example2} for results.

\begin{figure}[tbp]
\begin{center}
\begin{tabular}{cc}
\includegraphics[width=5cm]{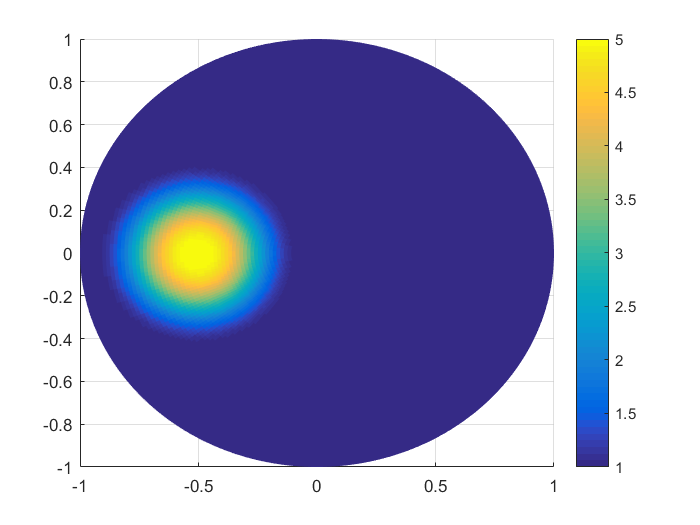} & %
\includegraphics[width=5cm]{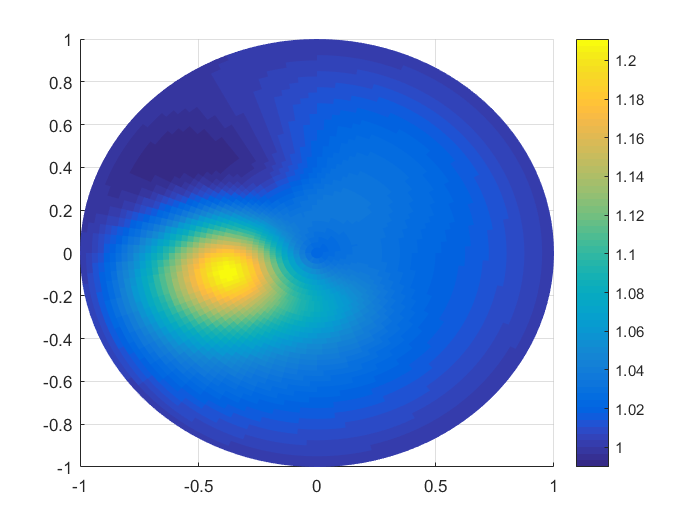} \\ 
(a) True $\sigma$ & (b) $N=4$ \\ 
\includegraphics[width=5cm]{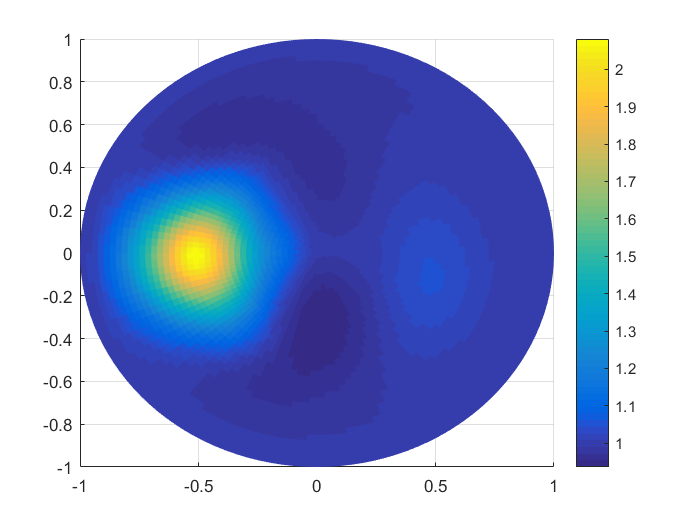} & %
\includegraphics[width=5cm]{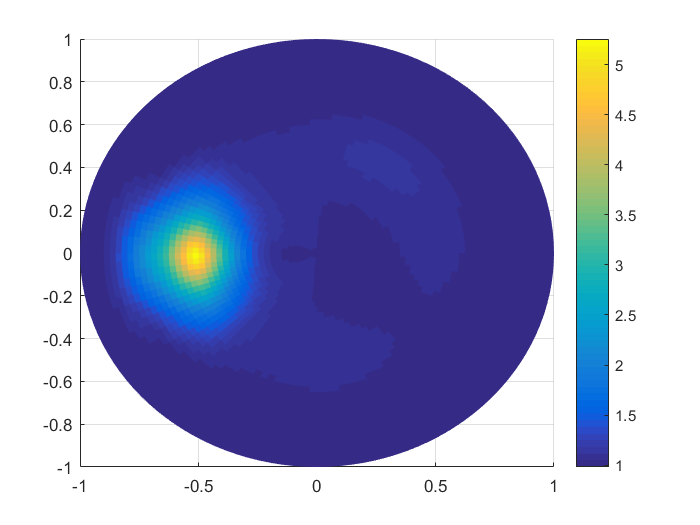} \\ 
(c) $N=6$ & (d) $N=8$%
\end{tabular}%
\end{center}
\caption{\emph{Results of Test 4. Imaging of a single inclusion with a high
inclusion/background contrast 5:1. Here, }$\protect\sigma =5$\emph{\ inside
the inclusion and }$\protect\sigma =1$ \emph{outside. a) Correct image. b)
Computed image for }$N=4$\emph{. c) Computed image for }$N=6$\emph{. d)
Computed image for }$N=8$\emph{, the best one out of three.}}
\label{example3}
\end{figure}

\textbf{Test 3}: In this example, we test the reconstruction method for a
single inclusion with a rather high inclusion/background contrast 5:1. The
results are shown on Figure \ref{example3}.

\begin{figure}[tbp]
\begin{center}
\begin{tabular}{cc}
\includegraphics[width=5cm]{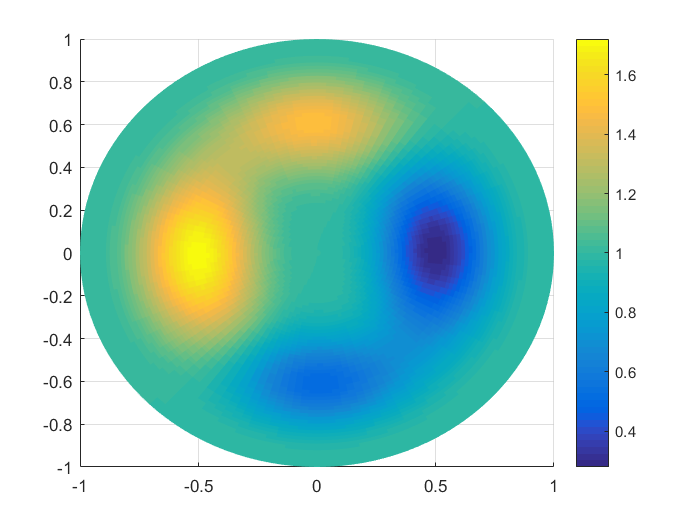} & %
\includegraphics[width=5cm]{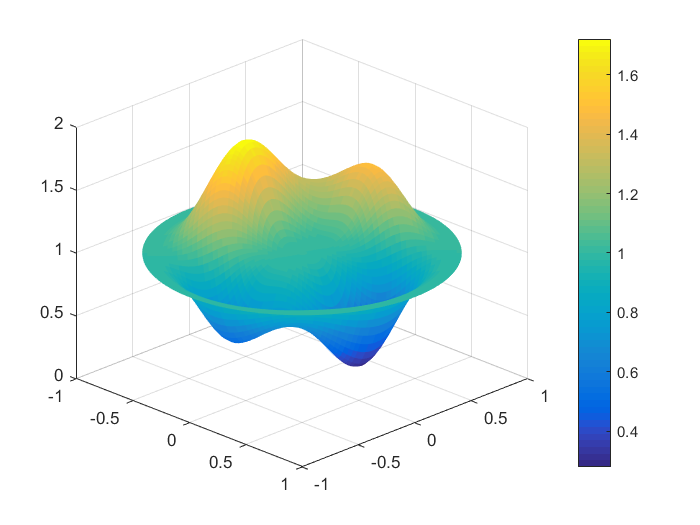} \\ 
(a) True $\sigma$ & (b) 3D view of true $\sigma$ \\ 
\includegraphics[width=5cm]{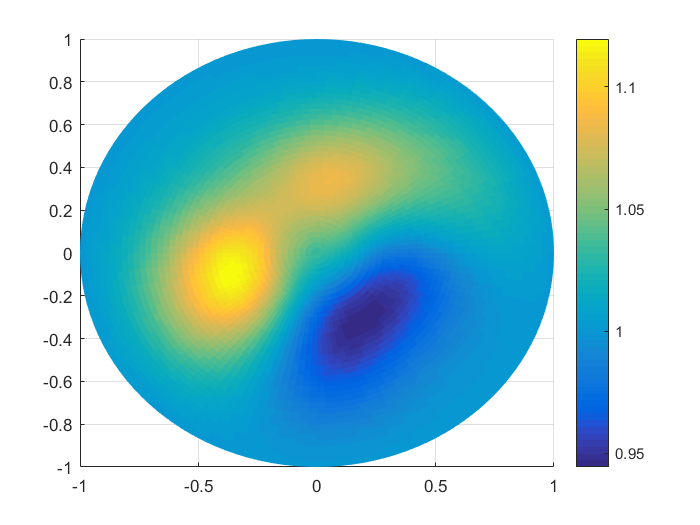} & %
\includegraphics[width=5cm]{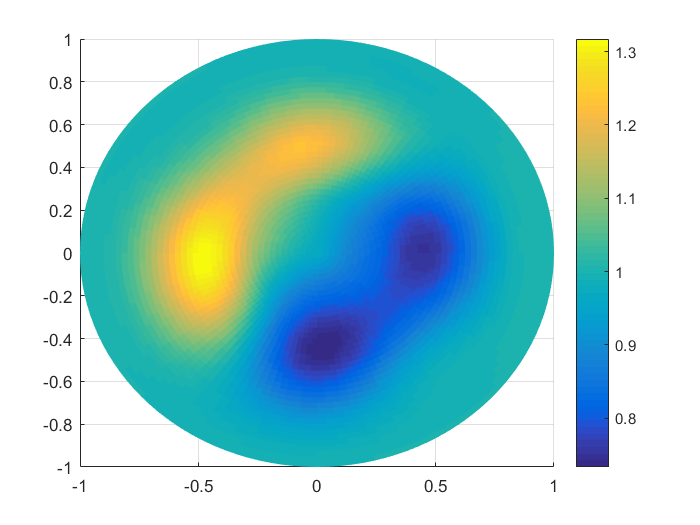} \\ 
(c) $N=4$ & (d) $N=6$ \\ 
\includegraphics[width=5cm]{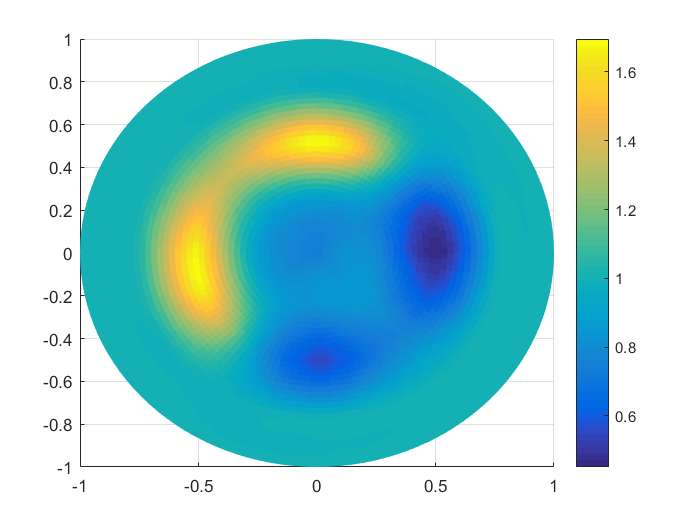} & %
\includegraphics[width=5cm]{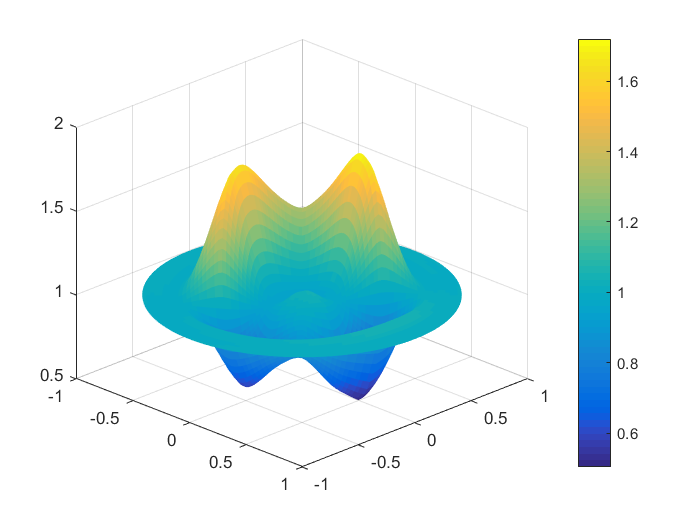} \\ 
(e) $N=8$ & (f) 3D view of the result $N=8$%
\end{tabular}%
\end{center}
\caption{\emph{Results of Test 4. We now test imaging of a smoothly varying
conductivity rather than of inclusions above. The values of }$\protect\sigma %
\left( x\right) $ \emph{inside of the inhomogenuity vary in a wide range }$%
\protect\sigma _{\min }\approx 0.3$ and $\protect\sigma _{\max }\approx 1.7.$
\emph{And }$\protect\sigma =1$\emph{\ in the homogeneous part of this disk. }
\emph{Here, we have stopped on the 3rd mesh refinement. a) Correct 2D image.
b) 3D presentation of a). c) Computed image for }$N=4$\emph{. d) Computed
image for }$N=6$\emph{. e) Computed 2D image for }$N=8$\emph{. f) 3D
presentation of e). Thus, we can accurately image not only \textquotedblleft
sharp" inclusions but smoothly varying functions as well. }}
\label{example4}
\end{figure}

\textbf{Test 4}. We now test our method for the case when the function $%
\sigma \left( x\right) $ is smoothly varying within an abnormality and with
a wide range of variations between 0.4 and 1.6. The results are shown in
Figure \ref{example4}. Again $N=8$ is the best value out of three listed in (%
\ref{6.4}). Thus, our method can accurately image not only \textquotedblleft
sharp" inclusions as in Tests 1-3, but smoothly varying functions as well.

\begin{figure}[tbp]
\begin{center}
\begin{tabular}{cc}
\includegraphics[width=5cm]{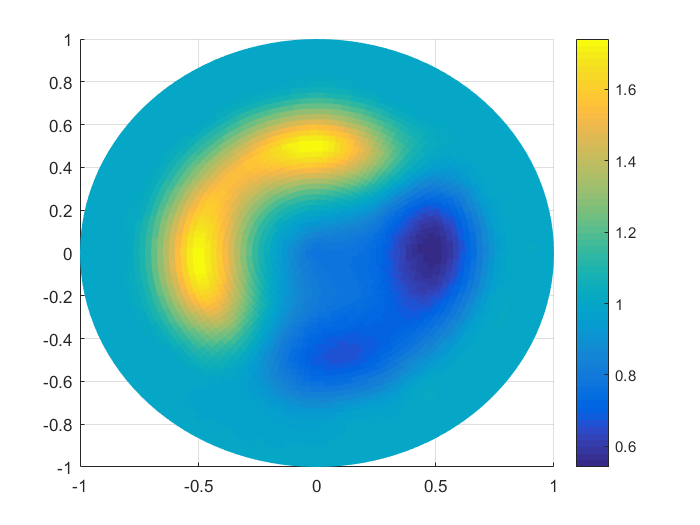} & %
\includegraphics[width=5cm]{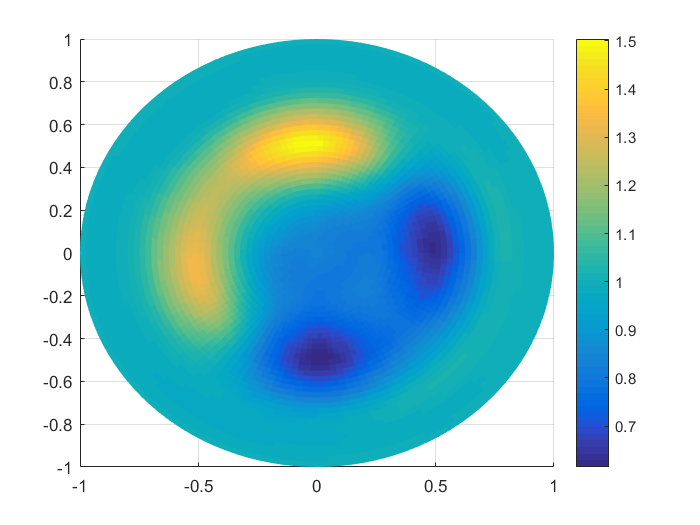} \\ 
(a) Noise level $1\%$ & (b) Noise level $10\%$.%
\end{tabular}%
\end{center}
\caption{\emph{Results of Test 5. In this test we have introduced random
noise in the data of test 4. Here, }$N=8$\emph{. a) Computed image with 1\%
noise. b) Computed image with 10\% noise.}}
\label{example5}
\end{figure}

\textbf{Test 5}. In this example we test the stability of the algorithm with
respect to the random noise in the data. We test the most challenging case
among ones above: the case of the function $\sigma \left( x\right) $ of Test
4. We set $N=8$. The noise is added for $x\in S_{1}$ and for the source $s$
as in (\ref{6.20}), $s\in \lbrack 0,2\pi ]$: 
\begin{equation}
g_{0,\text{noise}}(x,s)=g_{0}(x,s)(1+\epsilon \xi _{s})\text{ and }g_{1,%
\text{noise}}(x,s)=g_{1}(x,s)(1+\epsilon \xi _{s}),  \nonumber
\end{equation}%
where $\varepsilon $ is the noise level and $\xi _{s}$ is the independent
random variable depending only on the source position $s$ and uniformly
distributed on $[-1,1]$. The computational results are displayed on Figure %
\ref{example5} for the levels of noise of 1\% and 10\%.

This example indicates that our method is quite stable with respect to the
noise in the measured data.

\textbf{Test 6}. In all above tests 1-5 we have used the Dirichlet and
Neumann data on the entire boundary $S_{1}$ of our disk $\Omega .$ Also, the
source was running along the entire circle $C^{\left( s\right) }$ as in (\ref%
{6.20}). In this test, however, we study the case of incomplete data. First,
we work with the case when the source runs over the entire circle (\ref{6.20}%
) while the data $g_{0}\left( x,s\right) $ and $g_{0}\left( x,s\right) $ are
measured only on a part of the circle $S_{1}.$ Next, we study the case when
the source runs only along a part of the circle $C^{\left( s\right) }$ in (%
\ref{6.20}) while the data are measured on the entire circle $S_{1}.$ We
again use $N=8$ and the same function $\sigma \left( x\right) $ as in Test 4.

\begin{figure}[tbp]
\begin{center}
\begin{tabular}{cc}
\includegraphics[width=5cm]{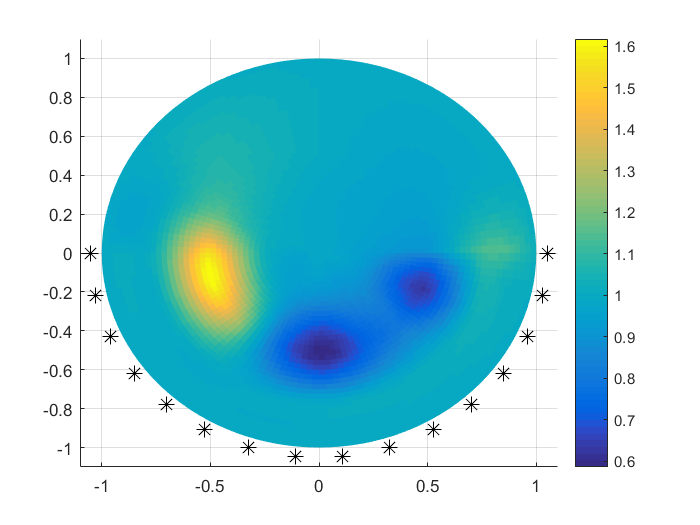} & %
\includegraphics[width=5cm]{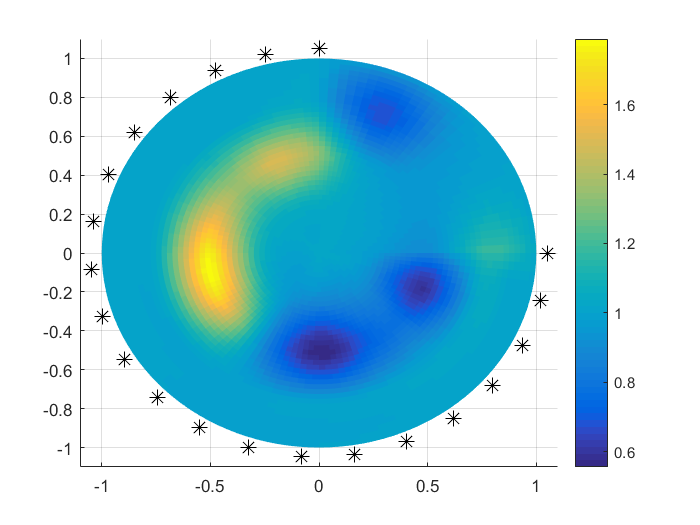} \\ 
(a) $50\%$ of the boundary data & (b) $75\%$ of the boundary data \\ 
\includegraphics[width=5cm]{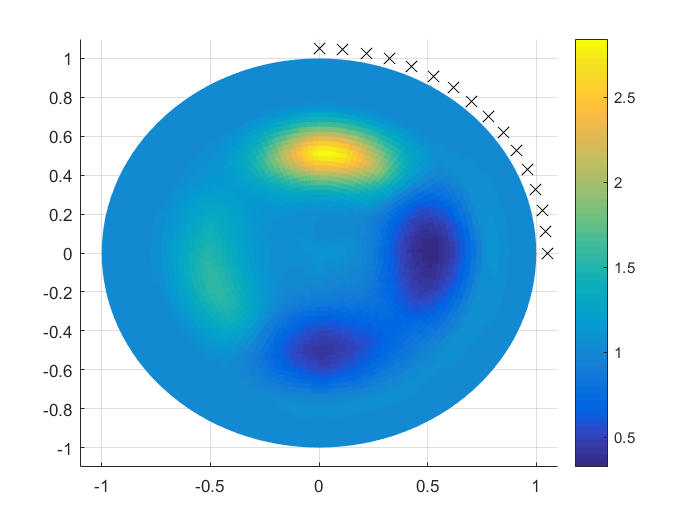} & %
\includegraphics[width=5cm]{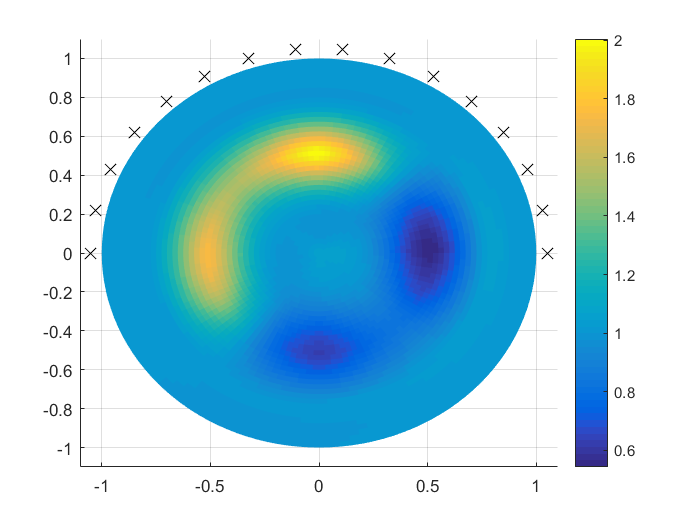} \\ 
(c) $25\%$ of the position of the source & (d) $50\%$ of the position of the
source%
\end{tabular}%
\end{center}
\caption{\emph{Results of Test 6. While the structure to be imaged is the
same as the one of Figure 4a), a lesser amount of data is used here. In a)
and b) we use incomplete boundary data on the circle }$S_{1}=\left\{
r=1\right\} $\emph{, while the source is still running as in (\protect\ref%
{6.20}): }$s\in \left( 0,2\protect\pi \right) ,$ i.e.\emph{\ over the entire
circle }$C^{\left( s\right) }=\left\{ r=4\right\} .$\emph{\ On the other
hand, in c) and d) the boundary data are measured at the entire circle }$%
S_{1}=\left\{ r=1\right\} ,$\emph{\ while the source is running over only a
part of the circle }$C^{\left( s\right) }=\left\{ r=4\right\} .$\emph{\ In
a) and b) }$\ast $\emph{\ indicates the part of the circle }$S_{1}=\left\{
r=1\right\} $\emph{\ where the data are measured. In c) and d) }$\times $%
\emph{\ indicates the part of the circle }$\left\{ r=4\right\} $\emph{\
where the source runs. The part with }$\times $\emph{\ is depicted on }$%
\left\{ r=1\right\} $\emph{\ rather than on }$\left\{ r=4\right\} $ \emph{%
only for the convenience of the presentation.}}
\label{example6}
\end{figure}

Figures \ref{example6} display results of Test 6. Comparing with the correct
image of Figure \ref{example4}, one can observe that, using 50\% of the
measured boundary data, one looses about 50\% of the internal information.
On the other hand, using 50\% of the positions of the source, one can still
recover the internal conductivity with a rather good accuracy. Hence, it
seems to be more important to measure at the entire boundary than to use the
entire circle $C^{\left( s\right) }$ for the positions of the source.

\section{Concluding remarks}

\label{sec:7}

Using a new concept, which was proposed in \cite{KlibJIIP}, we have
developed here the convexification numerical method for the inverse problem
of Electrical Impedance Tomography. While in all past publications on the
convexification \cite{BKconv,KlibSIMA97,KlibKam,KlibThanh,KlibKol1,KlibKol2}
only a single location of the source was used for either time dependent or
frequency dependent data, in the current paper the Dirichlet and Neumann
data are generated by a point source running along an interval of a straight
line and these data are independent on neither time nor frequency. We have
proved theorems, assuring the global convergence of our method. The key
analytical tool here is the tool of Carleman estimates. In particular, we
have proven two new Carleman estimates.

We have conducted extensive numerical testing of our method. Our
computational results demonstrate that this technique can accurately image
both sharp inclusions and smoothly varying abnormalities. In addition, our
method is quite stable with respect to the noise in the data (Test 5). We
have also studied numerically the performance of our method for the case
when either boundary data are measured only on a part of the boundary or the
source is running only on a part of the circle surrounding our domain of
interest.


\begin{thebibliography}{99}
\bibitem{Mueller2016} M.~Alsaker and J.~L. Mueller. \newblock {\em A d-bar
algorithm with a priori information for 2-dimensional electrical impedance
tomography,} \newblock SIAM Journal on Imaging Sciences,
9 (2016), pp. 1619-1654.

\bibitem{AmmariWenlong2016} H.~Ammari, G.~S. Alberti, B.~Jin, J.~K. Seo, and
W.~Zhang. \newblock {\em The linearized inverse problem in multifrequency
electrical impedance tomography,} 
\newblock SIAM Journal on Imaging
Sciences, 9 (2016), pp. 1525-1551.

\bibitem{anomaly1} H.~Ammari and H.~Kang. 
\newblock {\em Reconstruction of small inhomogeneities from boundary
  measurements}, Springer-Verlag, Berlin, 2004.

\bibitem{Wenlong2017} H.~Ammari, L.~Qiu, F.~Santosa, and W.~Zhang. \newblock 
{\em Determining anisotropic conductivity using diffusion tensor imaging data in
magneto-acoustic tomography with magnetic induction,} 
\newblock {Inverse
Problems}, 33 (2017), pp. 125006.

\bibitem{BakKlib} A.~B. Bakushinskii, M.~V. Klibanov and N.~A. Koshev, \emph{Carleman weight functions for a globally convergent numerical
method for ill-posed Cauchy problems for some quasilinear PDEs}, Nonlinear
Analysis: Real World Applications, 34 (2017), pp.~201--224.

\bibitem{Baud} L. Baudouin, M. de Buhan and S. Ervedoza, \emph{Convergent
algorithm based on Carleman estimates for the recovert of a potential in the
wave equation}, SIAM J. on Numerical Analysis, 55 (2017), pp. 1578-1613.

\bibitem{BK} \textsc{L.~Beilina and M.~V. Klibanov}, \emph{Approximate
Global Convergence and Adaptivity for Coefficient Inverse Problems},
Springer, New York, 2012.

\bibitem{BKconv} \textsc{L.~Beilina and M.~V. Klibanov}, \emph{Globally
strongly convex cost functional for a coefficient inverse problem},
Nonlinear Analysis: Real World Applications, 22  (2015), pp.~272--288.

\bibitem{BY} \textsc{M. Bellassoued} and \textsc{M. Yamamoto}, \emph{%
Carleman Estimates and Applications to Inverse Problems for Hyperbolic
Systems}, Springer Japan KK, 2017.

\bibitem{BukhKlib} \textsc{A.~Bukhgeim and M.~Klibanov}, \emph{Uniqueness in
the large of a class of multidimensional inverse problems}, Soviet Math.
Doklady, 17 (1981), pp.~244--247.

\bibitem{borcea} L.~Borcea. \newblock {\em Electrical impedance tomography,}
\newblock { Inverse Problems}, 18 (2002), pp. R99-R136.

\bibitem{Calderon:1980} A.~P. Calderon. \newblock {\em On an inverse boundary
value problem,}
\newblock {Seminar on Numerical Analysis and its Applications to Continuum
  Physics (Rio de Janeiro, 1980), Soc. Brasil. Mat., Rio de Janeiro}, pp.
65--73.

\bibitem{Mueller2014} M.~Dodd and J.~Mueller. \newblock {\em A real-time d-bar
algorithm for 2-d electrical impedance tomography data,}  Inverse Probl Imaging (Springfield).  8 (2014), pp. 1013-1031.

\bibitem{Gilbarg} D. Gilbarg and N.S. Trudinger, Elliptic Partial
Differential Equations of Second Order, Springer, New York, 1984.

\bibitem{Siltanen2014} S.~J. Hamilton, M.~Lassas, and S.~Siltanen. \newblock %
{\em A direct reconstruction method for anisotropic electrical impedance
tomography,} \newblock  Inverse Problems, 30 (2014), pp. 075007.

\bibitem{Siltanen2016} S.~J. Hamilton, J.~M. Reyes, S.~Siltanen, and
X.~Zhang. \newblock {\em A hybrid segmentation and d-bar method for electrical
impedance tomography,}  SIAM Journal on Imaging Sciences,
9 (2016), pp. 770--793.

\bibitem{Harrah} B. Harrah and M.N. Minh, \emph{Enhancing residual-based
techniques with shape reconstruction features in Electrical Impedance
Tomography,} Inverse Problems, 32 (2016), pp. 125002.

\bibitem{Holder:2005} D.~Holder. 
\newblock {\em Electrical {I}mpedance {T}omography: {M}ethods, {H}istory and
  {A}pplications}, Institute of Physics, Bristol, 2005.

\bibitem{Paivarinta2018} N. Hyv$\ddot{o}$nen , L. P$\ddot{a}$iv$\ddot{a}$%
rinta and J. P. Tamminen \newblock {\em Enhancing D-bar reconstructions for
electrical impedance tomography with conformal maps. Enhancing D-bar
reconstructions for electrical impedance tomography with conformal maps,}
Inverse Problems and Imaging, 12 (2018), pp. 373-400.

\bibitem{Bangtisparse2012} B.~Jin, T.~Khan, and P.~Maass. \newblock { \em A
reconstruction algorithm for electrical impedance tomography based on
sparsity regularization,} 
International Journal for Numerical
Methods in Engineering, 89 (2012), pp. 337--353.

\bibitem{BangtiFEM2016} B.~Jin, Y.~Xu, and J.~Zou. \newblock {\em A convergent
adaptive finite element method for electrical impedance tomography,} 
\newblock IMA Journal of Numerical Analysis, 37 (2017), pp. 1520--1550.

\bibitem{KlibSIMA97} \textsc{M.~V. Klibanov}, \emph{Global Convexity in a
Three-Dimensional Inverse Acoustic Problem}, SIAM Journal on Mathematical
Analysis, 28 (1997), pp.~1371--1388.

\bibitem{Ksurvey} \textsc{M.~V. Klibanov}, \emph{Carleman estimates for
global uniqueness, stability and numerical methods for coefficient inverse
problems}, Journal of Inverse and Ill-Posed Problems, 21 (2013),
pp.~477--560.

\bibitem{KQR} \textsc{M.~V. Klibanov}, \emph{Carleman estimates for the
regularization of ill-posed Cauchy problems, }Applied Numerical Mathematics,
94 (2015), pp. 46--74.

\bibitem{KlibCauchy} \textsc{M.~V. Klibanov}, \emph{Carleman weight
functions for solving ill-posed Cauchy problems for quasilinear PDEs},
Inverse Problems, 31 (2015), pp.~125007.

\bibitem{KlibThanh} \textsc{M.~V. Klibanov and N.~T. Th{\`{a}}nh}, \emph{%
Recovering Dielectric Constants of Explosives via a Globally Strictly Convex
Cost Functional}, SIAM Journal on Applied Mathematics, 75 (2015),
pp.~518--537.

\bibitem{KlibKam} \textsc{M.~V. Klibanov }and\textsc{\ V.~G. Kamburg}, \emph{%
Globally strictly convex cost functional for an inverse parabolic problem},
Mathematical Methods in the Applied Sciences, 39 (2016), pp.~930--940.

\bibitem{KlibYag} \textsc{M.~V. Klibanov, N.A. Koshev, J. Li }and\textsc{\
A.G. Yagola, }\emph{Numerical solution of an ill-posed Cauchy problem for a
quasilinear parabolic equation using a Carleman weight function}, J. Inverse
and Ill-Posed Problems, 24 (2016), pp. 761-776.

\bibitem{KlibKol1} \textsc{M.~V. Klibanov, A.~E. Kolesov, L.~Nguyen and
A.~Sullivan}, \emph{Globally strictly convex cost functional for a 1-D
inverse medium scattering problem with experimental data}, SIAM J. Applied
Mathematics, 77 (2017), pp. 1733-1755.

\bibitem{KlibKol2} \textsc{M.~V. Klibanov and A.~E. Kolesov, }\emph{%
Convexification of a 3-D coefficient inverse scattering problem}, accepted
in Computers and Mathematics with Applications, 2018, also see arxiv:
1801.04404.

\bibitem{KlibJIIP} \textsc{M.V. Klibanov}, \emph{Convexification of
restricted Dirichlet-to-Neumann map}, J. Inverse and Ill-Posed Problems, 25
(2017), pp. 669-685.

\bibitem{anomaly4} O.~Kwon, J.-K. Seo, and J.-R. Yoon. \newblock {\em A real-time
algorithm for the location search of discontinuous conductivities with one
measurement,} \newblock  Comm. Pure Appl. Math., 55 (2002), pp. 1--29.

\bibitem{Zouli2007} J.~Li and J.~Zou. \newblock {\em A multilevel model
correction method for parameter identification,} 
\newblock  Inverse
Problems, 23 (2007), pp. 1759.

\bibitem{LRS} M.M. Lavrentiev, V.G. Romanov and S.P. Shishatskii, \emph{%
Ill-Posed Problems of Mathematical Physics and Analysis}, AMS, Providence,
R.I., 1986.

\bibitem{BangtiFEM2014} G.~Matthias, J.~Bangti, and L.~Xiliang. \newblock {\em An
analysis of finite element approximation in electrical impedance tomography,} %
\newblock  Inverse Problems, 30 (2014), pp. 045013.

\bibitem{Nachman} \textsc{A. Nachman}, \emph{Global uniqueness for a
two-dimensional inverse boundary value problem}, Annals of Mathematics, 143
(1996), pp. 71-96.

\bibitem{Scales} \textsc{J.~A. Scales, M.~L. Smith, and T.~L. Fischer}, 
\emph{Global optimization methods for multimodal inverse problems}, Journal
of Computational Physics, 103 (1992), pp.~258--268.

\bibitem{SeoLeeKim:2008} J.~K. Seo, J.~Lee, S.~W. Kim, H.~Zribi, and E.~J.
Woo. 
\newblock {\em Frequency-difference electrical impedance tomography (fdEIT):
  algorithm development and feasibility study}, \newblock { Phys. Meas.},
29 (2008), pp. 929--944.

\bibitem{hybrid8} J.-K. Seo and E.-J. Woo. 
\newblock {\em Magnetic resonance
electrical impedance tomography (MREIT)}. \newblock { SIAM Rev.},
53 (2011), pp. 40--68.

\bibitem{SU} \textsc{J. Sylvester} and \textsc{G. Uhlmann}, \emph{Global
uniqueness theorem for an inverse boundary problem}, Annals of Mathematics,
39 (1987), pp. 91-112.

\bibitem{T} \textsc{A.N. Tikhonov, A.V. Goncharsky, V.V. Stepanov }and 
\textsc{\ A.G. Yagola}, \emph{Numerical Methods for the Solution of
Ill-Posed Problems}, Kluwer, London, 1995.

\bibitem{hybrid9} T.~Widlak and O.~Scherzer. \newblock {\em Hybrid tomography for
conductivity imaging,} \newblock  Inverse Problems, 28(2012), pp. 084008.
\end{thebibliography}
\end{document}